\documentclass[11pt, oneside]{article}
\usepackage{amsfonts}
\usepackage{mathrsfs}
\usepackage{amsfonts}
\usepackage{latexsym,amssymb,amsmath}
\usepackage{latexsym}
\usepackage{amssymb}
\usepackage{amsmath}
\usepackage{amsthm}
\usepackage{indentfirst}
\usepackage{color}
\usepackage{graphicx}
\usepackage{bm,enumerate}
\usepackage[square, comma, sort&compress, numbers]{natbib}
\usepackage[colorlinks]{hyperref}
\numberwithin{equation}{section} \allowdisplaybreaks

\makeatletter
\renewcommand\@cite[1]{#1\hspace{0.2em}}
\makeatother
\newtheorem{theorem}{\color{black}\indent Theorem}[section]
\newtheorem{lemma}[theorem]{\color{black}\indent Lemma}

\newtheorem{definition}[theorem]{\color{black}\indent Definition}
\newtheorem{remark}[theorem]{\color{black}\indent Remark}

\newtheorem{assumption}[theorem]{\color{black}\indent Assumption}

\usepackage{amssymb,amsmath}
\DeclareMathOperator{\dive}{div}
\def\rd{{\rm d}}
\def\re{{\rm e}}

\begin{document}	
\title{\Large\bf  Large-data solutions in multi-dimensional thermoviscoelasticity with temperature-dependent viscosities
	\author{$\rm{Chuang~Ma}$, $\rm{Bin~Guo}$
			\thanks{Corresponding author\newline \hspace*{6mm}{\it Email
				addresses:} machuang24@mails.jlu.edu.cn(Chuang Ma),~bguo@jlu.edu.cn~(Bin Guo)
			}
		\\[4pt]
		School of Mathematics, Jilin University, Changchun 130012, PR China
}}
\date{\today} \maketitle

\begin{abstract}
This paper investigates a quasilinear parabolic system arising in thermoviscoelasticity of Kelvin-Voigt type with temperature-dependent viscosity and coupled terms. The system, given by
\begin{equation*}
	\begin{cases}
		u_{tt}=\nabla\cdot\big(\gamma(\Theta)\nabla u_t\big)+a\Delta u-\nabla\cdot f(\Theta), & x \in \Omega,\ t > 0, \\
		\Theta_t=\Delta\Theta+\gamma(\Theta)|\nabla u_t|^2-f(\Theta)\nabla u_t, & x \in \Omega,\ t > 0, \\
		u=0,\quad\frac{\partial\Theta}{\partial\nu}=0, & x \in \partial\Omega,\ t > 0, \\
		u(x,0)=u_0(x),\; u_t(x,0)=u_{0t}(x),\;\Theta(x,0)=\Theta_0(x), & x \in \Omega,
	\end{cases}
\end{equation*}
models heat generation by acoustic waves in solid materials and can be derived as a scalar simplification of more complex piezoelectric-thermoviscoelastic model.
Under the assumptions that $u_0\in H_0^1(\Omega)$, $u_{0t}\in L^2(\Omega)$, $\Theta_0\in L^1(\Omega)$ with $\Theta_0\geqslant0$ a.e.~in $\Omega$, that $\gamma,f\in C^0([0,\infty))$ satisfy $f(0)=0$, and that there exist constants $k_\gamma,K_\gamma,K_f>0$ and $0<\alpha<\frac{N+2}{2N}$ such that
$$k_\gamma\leqslant\gamma(\xi)\leqslant K_\gamma\quad\text{and}\quad |f(\xi)|\leqslant K_f(1+\xi)^\alpha\qquad\forall~\xi\geqslant0,$$
we establish the global existence of weak solutions for arbitrarily large initial data in bounded domains $\Omega\subset\mathbb{R}^N$ ($N\geqslant1$). 
The result extends recent one-dimensional finding \cite{WinklerZAMP} to the multi-dimensional setting without requiring any smallness condition on the data.
\end{abstract}

\thispagestyle{empty}

\section{Introduction}
This paper investigates the following initial-boundary value problem
\begin{equation}\label{equ1}
	\begin{cases}
		u_{tt}=\nabla\cdot\big(\gamma(\Theta)\nabla u_t\big)+a\Delta u-\nabla\cdot f(\Theta) , & x \in \Omega,\ t > 0, \\
		\Theta_t=\Delta\Theta+\gamma(\Theta)|\nabla u_t|^2-f(\Theta)\nabla u_t, & x \in \Omega,\ t > 0, \\
		u=0,\quad\frac{\partial\Theta}{\partial\nu}=0, & x \in \partial\Omega,\ t > 0, \\
		u(x,0)=u_0(x),\; u_t(x,0)=u_{0t}(x),\;\Theta(x,0)=\Theta_0(x), & x \in \Omega,
	\end{cases}
\end{equation}
where $\Omega\subset\mathbb{R}^N$ ($N\geqslant 1$) is a bounded domain with smooth boundary $\partial\Omega$, and $a>0$ is a constant. 
This system \eqref{equ1} describes heat generation caused by acoustic waves in materials of Kelvin-Voigt type. 
The functions $u$ and $\Theta$ represent the displacement field and the temperature, respectively, while $\gamma$, $a$ and $f$ denote given coefficient functions whose precise properties will be specified later.
On one hand, when $\gamma$ is constant, \eqref{equ1} constitutes a classical model for the dissipation of mechanical energy into heat through viscous effects in solid materials \cite{RT2009,RZ1997}. On the other hand, it can be derived as a simplified version of a more elaborate model describing thermoviscoelastic processes in piezoelectric materials. A more general model arising in piezoelectric materials reads as follows:
\begin{align}\label{rhoutt}
	\begin{cases}
		\rho u_{tt}=\dive\left[\mathbb{C}(\nabla^s u-\Theta\mathbb{B})+\mathbb{D}(\nabla^s u_t)-e^T E\right],\\[2mm]
		c(\Theta)\Theta_t-\dive(\mathbb{K}\nabla\Theta)=\mathbb{D}(\nabla^s u_t):\nabla^s u_t-\Theta\mathbb{C}\mathbb{B}:\nabla^s u_t,
	\end{cases}
\end{align}
where the following notations in the model \eqref{rhoutt} are employed
\begin{align*}
&u: (x,t)\rightarrow\mathbb{R}^N~\text{displacement},\\
&\Theta: (x,t)\rightarrow\mathbb{R}~\text{temperature},\\
&\rho>0~\text{mass density},\\
&c(\Theta)>0~\text{heat capacity (depending on $\Theta$)},\\
&\nabla^s u=\frac12\big(\nabla u+(\nabla u)^{T}\big)~\text{the usual symmetric gradient},\\
&\mathbb{C}~\text{the 4th-order tensor of elastic constants},\\
&\mathbb{B}\in \mathbb{R}_{sym}^{N\times N}~\text{the thermal-dilation tensor},\\
&\mathbb{D}: \mathbb{R}_{sym}^{N\times N}\rightarrow \mathbb{R}_{sym}^{N\times N}~\text{a (possibly nonlinear) viscosity},\\
&\mathbb{K}\in  \mathbb{R}_{sym}^{N\times N}~\text{the heat-conductivity tensor},\\
&e~\text{the piezoelectric coupling tensor}, \\
&E~\text{the electric field strength}.
\end{align*}
This piezoelectric-thermoviscoelastic system is complemented by the constitutive relation for the electric displacement field
\begin{equation}\label{divD}
	D=\varepsilon E+e\nabla^s u,
\end{equation}
where $\varepsilon$ denotes the permittivity matrix. In the absence of free charges, $D$ is divergence-free, i.e. $\dive D=0$. Under the simplifying assumptions $\rho\equiv1$, $c(\Theta)\equiv1$, $\mathbb{K}=\mathbb{I}$, and the viscous response of the form $\mathbb{D}(\nabla^s u_t)=d(\Theta,\nabla^s u)\nabla^s u_t$ with the viscosity tensor $d$, we may eliminate the electric field $E$ by combining \eqref{divD} with the condition $\dive D=0$. Substituting the resulting expression into \eqref{rhoutt} yields
\begin{align}\label{utt1}
	\begin{cases}
		u_{tt}=\dive\left[a\nabla^s u+\gamma\nabla^s u_t-f\right],\\
		\Theta_t-\Delta\Theta=\gamma\nabla^s u_t:\nabla^s u_t-f:\nabla^s u_t,
	\end{cases}
\end{align}
where the effective coefficients are given by
$$a=\frac{\mathbb{C}}{\rho}+\frac{e^2}{\varepsilon\rho},
\quad\gamma=\frac{d}{\rho},
\quad\text{and}\quad 
f=\frac{\Theta\mathbb{C}\mathbb{B}}{\rho}.$$
The structure of \eqref{utt1} resembles that of \eqref{equ1}, but there is an essential difference: in \eqref{utt1} the displacement $u$ is a vector-valued function and its symmetric gradient $\nabla^s u$ is a matrix-valued field, whereas \eqref{equ1} involves only scalar quantities. Thus, \eqref{equ1} may be regarded as a scalar prototype that captures the principal analytical difficulties of the full tensorial system \eqref{utt1}.

The analysis of one-dimensional versions of such thermoviscoelastic systems has a long history, starting with the fundamental contributions of Dafermos \cite{Dafermos1982} and Dafermos--Hsiao \cite{DH1982}. For generalizations of the form
\begin{equation}\label{eqDgeneral}
	\begin{cases}
		u_{tt} = (\gamma(\Theta,u_x) u_{xt})_x + (a(\Theta,u_x) u_x)_x - (f(\Theta,u_x))_x, \\
		\Theta_t = \Theta_{xx} + \gamma(\Theta,u_x) u_{xt}^2 - f(\Theta,u_x) u_{xt},
	\end{cases}
\end{equation}
a rather complete picture has emerged when the coefficients $\gamma$ and $a$ are independent of temperature. For constant $\gamma$ and suitable classes of functions $a=a(u_x)$ and $f=\Theta F(u_x)$, the global existence of weak solutions was established in \cite{RZ1997}. Small-data global solvability results can be found in \cite{SY1995, ZS1989, Kim1983}, and under certain physically reasonable hypotheses even large-data smooth solutions have been constructed \cite{CH1994, GZ1999, SZZ1999}. Moreover, for such constant $\gamma$, results on long-time stabilization towards equilibrium are available \cite{HL1998, RZ1997}, and global solvability has also been proved when $\gamma$ depends on the strain $u_x$ \cite{JS1993, Watson2000}.

The situation becomes markedly more involved when the viscosity and elasticity coefficients depend on the temperature. Although experimental evidence indicates that both viscosity and elasticity may exhibit non-negligible variations with temperature \cite{Friesen2024, GBSD2007}, the corresponding mathematical theory remains far from complete. For systems of type \eqref{eqDgeneral} with $\gamma=\gamma(\Theta)$, $a=a(\Theta)$ and $f=f(\Theta)$, only local existence for arbitrary data is available, while global solutions have been obtained exclusively only in settings of suitably small data \cite{WinklerAMO, Fricke2025, CWPrepint}.

Turning to higher-dimensional case, the situation is even less satisfactory. Even for the vectorial model \eqref{utt1} with temperature-independent coefficients, no global existence result in higher-dimensional case for large data seems to exist in the literature. When temperature-dependent coefficients are admitted, the difficulties multiply. Roub\'i\v{c}ek \cite{RT2009} proved the existence of very weak solutions for a $d$-dimensional thermoviscoelastic system of Kelvin-Voigt type featuring a nonlinear monotone viscosity of $p$-Laplacian form and a temperature-dependent heat capacity, subject to the structural conditions $p\geqslant2$, $\omega\geqslant1$, and $p>1+d/(2\omega)$. While this result accommodates arbitrarily large data with $L^1$-type regularity, it relies on specific growth restrictions on the heat capacity.
Even for the simpler situation where $\gamma$ does not depend on $\Theta$, the higher-dimensional system \eqref{utt1} poses substantial challenges at the level of basic solvability theory. Existing results either impose strong restrictions on the growth of $f$ (typically sublinear), require smallness of the data, or incorporate additional dissipative mechanisms \cite{BDGO, GJZW, MART, PIZW, RT2009}. For the special case $f\equiv0$, a recent contribution \cite{WinklerMMM} establishes global generalized solvability in higher-dimensional case. 

In the one-dimensional setting with temperature-dependent viscosities, Winkler \cite{WinklerZAMP} recently proved global existence of weak solutions for arbitrarily large initial data under the condition that the coupling function $f$ grows at most with exponent $\alpha<3/2$.
The purpose of the present paper is to extend the one-dimensional result of Winkler \cite{WinklerZAMP} to the $N$-dimensional case. More precisely, we establish global existence of weak solutions to \eqref{equ1} for arbitrary initial data belonging to suitable Sobolev space, under the assumptions that
$$\gamma=\gamma(\Theta),\quad a= const,\quad\text{and}\quad f=f(\Theta),$$
without any small data condition.
We shall approach this problem in the context of weak solvability, a concept that is made precise in the following definition.
\begin{definition}\label{solution}
We say that $(u,\Theta)$ is a weak solution to \eqref{equ1} if the following hold:
\begin{enumerate}[$\bullet$]
	\item The initial data is of regularity
	$$u_0\in H_0^1(\Omega),\quad u_{0t}\in L^2(\Omega),\quad \Theta_0\in L^1(\Omega),$$
	with $\Theta_0\geqslant0$ for all $x\in\Omega$.
	\item Solutions $\Theta$ and $u$ satisfy $\Theta \geqslant 0$ in $\Omega \times (0,\infty)$ and
	\begin{equation}\label{regular}
	\begin{cases}
		u\in C\left([0,T],L^2(\Omega)\right)\cap L^\infty\left([0,\infty); W_0^{1,2}(\Omega)\right),\\
		u_t\in L^\infty\left((0,\infty); L^2(\Omega)\right)\cap	L_{loc}^2\left([0,\infty); W_0^{1,2}(\Omega)\right)\cap L_{loc}^2\left([0,\infty)\times\Omega\right),\\
		\Theta\in\bigcap_{q\in\big[1,\frac{N+2}{N}\big)} L_{loc}^q\left([0,\infty)\times\Omega\right)\cap \bigcap_{r\in\big(1,\frac{N+2}{N+1}\big)}L_{loc}^r\left([0,\infty); W^{1,r}(\Omega)\right).
	\end{cases}
	\end{equation}
	\item The momentum equation
	\begin{align}\label{equmoment}
		-\int_{0}^{\infty} \int_{\Omega} &u_t \varphi_t\rd x\rd t
		- \int_{\Omega} u_{0t} \varphi(\cdot, 0)\rd x\\
		&=-\int_{0}^{\infty} \int_{\Omega} \gamma(\Theta)\nabla u_t \nabla\varphi \rd x\rd t
		-a \int_{0}^{\infty} \int_{\Omega} \nabla u \nabla\varphi \rd x\rd t
		+ \int_{0}^{\infty} \int_{\Omega} f(\Theta) \nabla\varphi\rd x\rd t,\notag
	\end{align}
	is satisfied for all $\varphi\in C_0^\infty(\Omega\times[0,\infty))$.	
	\item The temperature equation
	\begin{align}\label{equentropy}
		-\int_0^\infty\int_\Omega&\Theta\varphi_t\rd x\rd t-\int_0^\infty\int_\Omega\Theta_0\varphi(\cdot,0)\rd x\rd t\\
		&=-\int_0^\infty\int_\Omega\nabla\Theta\cdot\nabla\varphi\rd x\rd t
		+\int_0^\infty\int_\Omega\gamma(\Theta)|\nabla v|^2\varphi\rd x\rd t-\int_0^\infty\int_\Omega f(\Theta)\nabla u_t\varphi\rd x\rd t.\notag
	\end{align}
	holds for every $\varphi\in C_0^\infty(\Omega\times[0,\infty))$.
\end{enumerate}
\end{definition}

Next, we introduce the main assumptions that will be used throughout the paper.
\begin{assumption}\label{assum1}
	In this paper, we always assume that the initial data satisfy
	$$u_0\in H_0^1(\Omega),\quad u_{0t}\in L^2(\Omega),\quad\Theta_0\in L^1(\Omega),$$ 
	with $\Theta_0\geqslant0$ for all $x\in\Omega$.
	Moreover, the functions $\gamma\in C^0([0,\infty))$ and $f\in C^0([0,\infty))$ satisfy $f(0)=0$,
	 and 
	$$k_\gamma\leqslant \gamma(\xi)\leqslant K_\gamma,\quad\forall~\xi\geqslant0,$$
	and 
	$$|f(\xi)|\leqslant K_f(1+\xi)^\alpha,\quad\forall~\xi\geqslant0,$$
	with constants $k_\gamma>0,~K_\gamma>0,~K_f>0$ and $0<\alpha<\frac{N+2}{2N}$.
\end{assumption}

Our main result is the following theorem.
\begin{theorem}\label{Th1}
Let Assumption \ref{assum1} hold. Then there exists a pair $(u,\Theta)$ such that $\Theta \geqslant 0$ a.e. in $\Omega \times (0,\infty)$ and $(u,\Theta)$ is a global weak solution of \eqref{equ1} in the sense of Definition \ref{solution}.
\end{theorem}

\textbf{Outline of the paper.} 
The remainder of this paper is organized as follows.

In Section~2, we introduce a parabolic regularization of the original system \eqref{equ1} by adding fourth-order dissipation terms. This regularized problem \eqref{equ1*} is shown to be locally well-posed within the framework of classical solutions from  Lemma \ref{Lemlocal}, and a basic energy identity is derived (Lemma \ref{Lem-energy*}), which provides uniform estimates independent of the regularization parameter $\varepsilon$.

Section~3 is devoted to proving that for each fixed $\varepsilon\in(0,1)$, the regularized problem admits a global solution. By employing fractional power estimates and semigroup technique, we establish higher-order regularity properties that exclude the possibility of finite-time blow-up, leading to the global existence result stated in Lemma \ref{Lemglobal}.

In Section~4, we derive a series of $\varepsilon$-independent estimates that are essential for the limiting procedure. These include interpolation inequalities for the temperature (Lemma \ref{Lemq}), a Boccardo-Gallou\"et type estimate for the temperature gradient (Lemma~\ref{LemThetanabla}), and estimates for the time derivatives in dual spaces (Lemma~\ref{Lem*vuTheta}). Based on these uniform bounds, we extract convergent subsequences and identify the limit functions $(v,u,\Theta)$ in Lemma~\ref{Lemwekav}, along with their convergence properties and the weak formulation satisfied by $v$ and $u$.

Section~5 focuses on the solution properties of the temperature $\Theta$. We first introduce the Steklov averaging technique (Lemma~\ref{LemSteklov}) to handle time derivatives in the weak formulation. Then, in Lemma~\ref{Lemgammato}, we establish the strong convergence of  $\sqrt{\gamma_\varepsilon(\Theta_\varepsilon)}\nabla v_\varepsilon$ in $L^2$, which is crucial for passing to the limit in the quadratic nonlinearity $\gamma(\Theta)|\nabla v|^2$. Finally, combining all the convergence results, we complete the proof of Theorem~\ref{Th1}, thereby establishing the global existence of weak solutions to the original problem \ref{equ1}.

\section{Preliminaries}  
Let $\| \cdot \|_{p}$  denote the standard norm of $L^{p}(\Omega)$, and let $C$ represent a positive constant that may change at each occurrence.
To construct solutions of the original problem \eqref{equ1}, we will employ an approximation scheme based on a parabolic regularization. This approach is motivated by the observation that the substitution $v:=u_t$ formally transforms \eqref{equ1} into the following system
\begin{equation}\label{equvt}
	\begin{cases}
		v_{t}=\nabla\cdot\big(\gamma(\Theta)\nabla v\big)+a\Delta u-\nabla\cdot f(\Theta),~~& x \in \Omega,\ t > 0, \\
		u_t=v,~~& x \in \Omega,\ t > 0, \\
		\Theta_t=\Delta\Theta+\gamma(\Theta)|\nabla v|^2-f(\Theta)\nabla v,& x \in \Omega,\ t > 0, \\
		v=0,\quad u=0,\quad\frac{\partial\Theta}{\partial\nu}=0,& x \in \partial\Omega,\ t > 0, \\
		v(x,0)=u_{0t}(x),~u(x,0)=u_0(x),~\Theta(x,0)=\Theta_0(x),& x \in \Omega,
	\end{cases}
\end{equation}
which resembles a quasilinear parabolic system for lacking sufficient parabolicity in the following equation for $(v,u,\Theta)$ to apply standard theory directly.

To overcome this, we introduce a family of semilinear parabolic problems, parameterized by $\varepsilon\in(0,1)$, which are designed to be globally solvable and to preserve the essential energy structure of the original system \eqref{equvt}.  Specifically, we introduce a artificial fourth-order dissipation term $-\varepsilon\Delta^2 v$ in the equation for $v$ and a second-order dissipation term $\varepsilon\Delta u$ in the equation for $u$. These terms provide the necessary higher-order regularization to establish global existence for the corresponding approximate problems while vanishing in the limit $\varepsilon\searrow0$.

We first construct smooth approximations of the initial data and the nonlinearities $\gamma$ and $f$. Let $\{v_{0\varepsilon}\}_{\varepsilon\in(0,1)}\subset C_0^\infty([0,\infty))$, $\{u_{0\varepsilon}\}_{\varepsilon\in(0,1)}\subset C_0^\infty([0,\infty))$ and $\{\Theta_{0\varepsilon}\}_{\varepsilon\in(0,1)}\subset C^\infty([0,\infty))$ be such that $\Theta_{0\varepsilon}\geqslant0$ in $\Omega$ for all $\varepsilon\in(0,1)$, and that 
\begin{align}\label{smoothdata}
	\begin{cases}
		v_{0\varepsilon}\to v_0\quad\text{in}~ L^2(\Omega),\\
		u_{0\varepsilon}\to u_0\quad\text{in}~ W^{1.2}(\Omega),\\
		\Theta_{0\varepsilon}\to \Theta_0\quad\text{in}~ L^1(\Omega),
	\end{cases}
\end{align}
as $\varepsilon\searrow0$. Furthermore, let $\{\gamma_\varepsilon\}_{\varepsilon\in(0,1)}\subset C^\infty([0,\infty))$ and $\{f_\varepsilon\}_{\varepsilon\in(0,1)}\subset C^\infty([0,\infty))\cap L^\infty([0,\infty))$ be sequences approximating $\gamma$ anf $f$, respectively, such that 
\begin{equation}\label{gamma*bound}
	k_\gamma\leqslant\gamma_\varepsilon(\xi)\leqslant K_\gamma,
	\quad\forall~\xi\geqslant0,~\forall~\varepsilon\in(0,1),
\end{equation}
that
\begin{equation}\label{falpha}
	|f_\varepsilon(\xi)|\leqslant K_f(1+\xi)^\alpha,
	\quad\forall~\xi\geqslant0,~\forall~\varepsilon\in(0,1),
\end{equation}
with $k_\gamma>0,~K_\gamma>0,~K_f>0$ and $0<\alpha<\frac{N+2}{2N}$, and that
$$\gamma_\varepsilon\to \gamma\quad\text{and}\quad f_\varepsilon\to f\quad\text{in}~L_{loc}^\infty([0,\infty)),$$
as $\varepsilon\searrow0$. 
For each $\varepsilon\in(0,1)$, we then consider the following regularized parabolic problem
\begin{equation}\label{equ1*}
	\begin{cases}
		v_{\varepsilon t}= -\varepsilon\Delta^2v_{\varepsilon}+\nabla\cdot\big(\gamma_\varepsilon(\Theta_\varepsilon)\nabla v_{\varepsilon}\big)+a \Delta u_{\varepsilon}-\nabla\cdot f_\varepsilon(\Theta_\varepsilon), & x \in \Omega,\ t > 0, \\
		u_{\varepsilon t}=\varepsilon\Delta u_{\varepsilon}+v_\varepsilon, & x \in \Omega,\ t > 0, \\
		\Theta_{\varepsilon t}=\Delta\Theta_{\varepsilon}+\gamma_\varepsilon(\Theta_\varepsilon)|\nabla v_{\varepsilon}|^2-f_\varepsilon(\Theta_\varepsilon)\nabla v_{\varepsilon},& x \in \Omega,\ t > 0, \\
		v_\varepsilon=\Delta v_{\varepsilon}=0,~u_\varepsilon=0,~\frac{\partial\Theta_{\varepsilon}}{\partial\nu}=0, & x \in \partial\Omega,\ t > 0,\\
		v_\varepsilon(x,0)=v_{0\varepsilon}(x),~u_\varepsilon(x,0)=u_{0\varepsilon}(x),~\Theta_\varepsilon(x,0)=\Theta_{0\varepsilon}(x),~&x\in\Omega.
	\end{cases}
\end{equation}

Indeed, the regularized parabolic problem \eqref{equ1*} fits within the scope of the standard local solvability theory for parabolic systems (cf. \cite{Amann1993}), an approach similarly employed in \cite{WinklerMMM}, 
which guarantees a local-in-time classical solution in the following sense.
\begin{lemma}\label{Lemlocal}
Let $\varepsilon\in(0,1)$. Then there exist $T_{\max,\varepsilon}\in(0,\infty]$ and uniquely determined functions
	\begin{align*}
		\begin{cases}
			v_\varepsilon\in C^0\left(\overline{\Omega}\times[0,T_{\max,\varepsilon})\right)\cap C^{4,1}\left(\overline{\Omega}\times(0,T_{\max,\varepsilon}),\right),\\
			u_\varepsilon\in C^0\left(\overline{\Omega}\times[0,T_{\max,\varepsilon})\right)\cap C^{2,1}\left(\overline{\Omega}\times(0,T_{\max,\varepsilon}),\right)\cap C^0\left([0,T_{\max,\varepsilon}); W_0^{1,2}(\Omega)\right),\\
			\Theta_\varepsilon\in C^0\left(\overline{\Omega}\times[0,T_{\max,\varepsilon})\right)\cap C^{2,1}\left(\overline{\Omega}\times(0,T_{\max,\varepsilon}),\right),
		\end{cases}
	\end{align*}
such that $\Theta_\varepsilon\geqslant0$ in $\overline{\Omega}\times[0,T_{\max,\varepsilon})$, and \eqref{equ1*} is solved in the classical sense in  $\Omega\times(0,T_{\max,\varepsilon})$, and that
if $T_{\max,\varepsilon}<\infty$, then
	\begin{align}\label{Tmax}
		\limsup\limits_{t\rightarrow T_{\max,\varepsilon}}\left\{\|u_\varepsilon\|_{W^{1+\eta,\infty}(\Omega)}+\|\Theta_\varepsilon\|_{W^{1+\eta,\infty}(\Omega)}+\|v_\varepsilon\|_{W^{2+2\eta,\infty}(\Omega)}\right\}=\infty,~~~\forall~\eta>0.
	\end{align}
\end{lemma}
We begin by testing the first equation in \eqref{equ1*} with $v_{\varepsilon}$. After integrating the resulting expression and combining it with the remaining equations of the system, we obtain the following energy structure.
\begin{lemma}\label{Lem-energy*}
Regular solution $(v_\varepsilon, u_\varepsilon, \Theta_\varepsilon)$ of \eqref{equ1*} satisfies
\begin{align*} 
	\frac{\rd}{\rd t}\left\{\frac12\int_\Omega v_\varepsilon^2\rd x+\frac a2\int_\Omega|\nabla u_\varepsilon|^2\rd x+\int_\Omega\Theta_\varepsilon\rd x\right\}+\varepsilon\int_\Omega|\Delta v_{\varepsilon}|^2\rd x+\varepsilon a\int_\Omega|\Delta u|^2\rd x=0, 
\end{align*}
and $$\frac12\int_\Omega v_\varepsilon^2\rd x+\frac a2\int_\Omega|\nabla u_\varepsilon|^2\rd x+\int_\Omega\Theta_\varepsilon\rd x \leqslant \frac12\int_\Omega v_{0\varepsilon}^2\rd x+\frac a2\int_\Omega|\nabla u_{0\varepsilon}|^2\rd x+\int_\Omega\Theta_{0\varepsilon}\rd x.$$
\end{lemma}
\begin{proof}
Multiplying the first equation of \eqref{equ1*} by $v_\varepsilon$ and integrating over $\Omega$, it yields
\begin{align}\label{energy*v}
	\frac{d}{dt} \bigg\{ \frac{1}{2} \int_\Omega v_{\varepsilon}^2\rd x + \frac{a}{2} \int_\Omega |\nabla u_{\varepsilon}|^2\rd x\bigg\}&+\varepsilon\int_\Omega |\Delta v_{\varepsilon}|^2\rd x+\varepsilon a\int_\Omega |\Delta u_{\varepsilon}|^2\rd x\notag\\
	&=-\int_{\Omega}\gamma_\varepsilon(\Theta_\varepsilon)|\nabla v_{\varepsilon}|^2\rd x+\int_{\Omega}f_\varepsilon(\Theta_\varepsilon)\nabla v_{\varepsilon}\rd x.
\end{align}
Integrating the third equation of \eqref{equ1*} over $\Omega$ directly, we get
\begin{align*}
	\frac{\rd}{\rd t}\int_\Omega\Theta_{\varepsilon}\rd x =\int_\Omega\gamma_{\varepsilon}(\Theta_\varepsilon)|\nabla v_{\varepsilon}|^2\rd x-\int_{\Omega}f_\varepsilon(\Theta_\varepsilon)\nabla v_{\varepsilon}\rd x.
\end{align*}
Add the above two inequalities together. Notice that  the terms on the right hand side cancel. Hence we have \begin{align*} 
	\frac{\rd}{\rd t}\left\{\frac12\int_\Omega v_\varepsilon^2\rd x+\frac a2\int_\Omega|\nabla u_\varepsilon|^2\rd x+\int_\Omega\Theta_\varepsilon\rd x\right\}+\varepsilon\int_\Omega|\Delta v_{\varepsilon}|^2\rd x+\varepsilon a\int_\Omega|\Delta u|^2\rd x=0.
\end{align*}
Further, integrating the above formula over $(0,t)$, together with \eqref{smoothdata}, we get
\begin{align*}
	\frac12\int_\Omega v_\varepsilon^2\rd x+\frac a2\int_\Omega|\nabla u_\varepsilon|^2\rd x&+\int_\Omega\Theta_\varepsilon\rd x+\varepsilon\int_0^t\int_\Omega|\Delta v_{\varepsilon}|^2\rd x\rd t+\varepsilon a\int_0^t\int_\Omega|\Delta u|^2\rd x\rd t\\
	&=\frac12\int_\Omega v_{0\varepsilon}^2\rd x+\frac a2\int_\Omega|\nabla u_{0\varepsilon}|^2\rd x+\int_\Omega\Theta_{0\varepsilon}\rd x.
\end{align*}
Obviously, Lemma \ref{Lem-energy*} as follows by dropping some negative terms in the above formula.
\end{proof}

\section{Global solvability of the regularized parabolic problem}
The goal of this section is to establish the global solvability of the regularized problem \eqref{equ1*} for each fixed \(\varepsilon\in(0,1)\), i.e., to prove that \(T_{\max}=\infty\). 
To achieve this, it is sufficient to rule out any finite-time blow-up.
Our approach focuses on deriving higher-order regularity properties (which may depend on $\varepsilon$) by exploiting two key features: the artificial fourth-order diffusion mechanism in \eqref{equ1*} and the fact that $(f_\varepsilon)_{\varepsilon\in(0,1)}\subset L^\infty([0,\infty))$. These estimates will serve as the foundation for extending the solution globally in time.
A key tool for this purpose is a result from \cite{Winkler2005JDE}, which provides a fundamental estimate for the quantities at hand.

\begin{lemma}\label{LemA1}
For $p\geqslant 2$, let $A_1$ be the realization of $\varepsilon\Delta^2$ in $L^p(\Omega; \mathbb{R}^n)$ with domain $D(A_1)=\left\{ \psi \in W^{4,p}(\Omega; \mathbb{R}^n) \mid \psi = \Delta\psi = 0 \text{ on } \partial\Omega\right\}$.
Suppose $T_{\text{max}, \varepsilon} < \infty$ for some $\varepsilon \in (0,1)$. Then, for any $\beta \in (\frac12, \frac34)$, there exists a constant $C(\sigma, \varepsilon) > 0$ such that the corresponding fractional power satisfies
$$\|A_1^\beta v_\varepsilon(\cdot, t)\|_{L^p(\Omega)}\leqslant C(\sigma, \varepsilon),~\quad\forall~t \in (0, T_{\max,\varepsilon}).$$
\end{lemma}
\begin{proof}
Recalling the first equation in \eqref{equ1*}, we get
\begin{align*}
	v_{\varepsilon t}+A_1v_{\varepsilon}=\nabla\cdot\big(\gamma_\varepsilon(\Theta_\varepsilon)\nabla v_{\varepsilon}+a \nabla u_{\varepsilon}-f_\varepsilon(\Theta_\varepsilon)\big).
\end{align*}
With the help of a duality argument of Lemma 2.1 in \cite{Winkler2005JDE}, there exists $\beta^*=\beta^*(\beta) \in (\frac14, \frac34)$ such that
\begin{align}\label{vA1}
\|A_1^\beta v_\varepsilon\|_p
&=\left\|\re^{-tA_1}A^\beta v_{0\varepsilon}+\int_0^tA_1^\beta\re^{-(t-s)A_1}h_{\varepsilon}\rd s\right\|_p\notag\\
&\leqslant\|\re^{-tA_1}A_1^\beta v_{0\varepsilon}\|_p+\int_0^t\left\|A_1^\beta\re^{-(t-s)A_1}\dive\left[\gamma_\varepsilon(\Theta_\varepsilon)\nabla v_{\varepsilon}+ a\nabla u_{\varepsilon}-f_\varepsilon(\Theta_\varepsilon)\right]\right\|_p\rd s\notag\\
&\leqslant C\|A_1^\beta v_{0\varepsilon}\|_p+C_\beta\int_0^t(t-s)^{-\beta^*-\frac14}
\|\gamma_\varepsilon(\Theta_\varepsilon)\nabla v_{\varepsilon}+ a\nabla u_{\varepsilon}-f_\varepsilon(\Theta_\varepsilon)\|_p\rd s,
\end{align}
where
\begin{align*}
	\|\gamma_\varepsilon(\Theta_\varepsilon)\nabla v_{\varepsilon}+ a\nabla u_{\varepsilon}-f_\varepsilon(\Theta_\varepsilon)\|_p
	\leqslant C\|v_{\varepsilon}\|_{W^{1,p}(\Omega)}+C\|u_\varepsilon\|_{W^{1,p}(\Omega)}+C\|f_\varepsilon(\Theta_\varepsilon)\|_{L^\infty(\Omega)}.
\end{align*}	
We shall estimate the right-hand side of the above equality. The Gagliardo-Nirenberg inequality gives us that
\begin{align*}
	\|v_{\varepsilon}\|_{W^{1,p}(\Omega)} 
	&\leqslant C\|v_{\varepsilon}\|_{W^{4\beta,p}(\Omega)}^{\delta_1}
	\cdot\|v_{\varepsilon}\|_{L^p(\Omega)}^{1-\delta_1}\\
	&\leqslant C\|v_{\varepsilon}\|_{W^{4\beta,p}(\Omega)}^{\delta_1}
	\cdot\left(\|v_{\varepsilon}\|_{W^{1,p}(\Omega)}^{\delta_2}\cdot\|v_{\varepsilon}\|_{L^2(\Omega)}^{1-\delta_2}\right)^{1-\delta_1}\\
	&\leqslant C\|v_{\varepsilon}\|_{W^{4\beta,p}(\Omega)}^{\delta_1}\cdot \|v_{\varepsilon}\|_{W^{1,p}(\Omega)}^{(1-\delta_1)\delta_2},
\end{align*}
where $\delta_1=\frac{1}{4\beta},~\delta_2=\frac{N(p-2)}{2p+N(p-2)}$. 
We set $$M_{1\varepsilon}(T):=\sup\limits_{t\in(0,T)}\|A_1^\beta v_\varepsilon(\cdot,t)\|_p,~\quad\text{for}~T\in(0,T_{\max,\varepsilon}).$$
Thus, from the above inequality, we have 
\begin{align}\label{v*}
	\|v_{\varepsilon}(\cdot,t)\|_{W^{1,p}(\Omega)} \leqslant C\|A_1^\beta v_\varepsilon(\cdot,t)\|_{L^p(\Omega)}^{\delta_3}+C\leqslant M_{1\varepsilon}^{\delta_3}(T)+C,
\end{align}
where $\delta_3=\frac{\delta_1}{1-(1-\delta_1)\delta_2}<1$. 
Let $A_2$ be the realization of $-\varepsilon\Delta$ in $L^p(\Omega; \mathbb{R}^n)$ with domain $D(A_2)=\left\{ \psi \in W^{2,p}(\Omega; \mathbb{R}^n) \mid \psi = 0 \text{ on } \partial\Omega\right\}$. By using the second equation in the approximate problem \eqref{equ1*}, \eqref{v*} and the fact that $u_{0\varepsilon}\in C_0^\infty(\Omega)\subset D(A_2)$, we obtain that
\begin{align}\label{u*}
	\|u_\varepsilon\|_{W^{1,p}(\Omega)}
	&=\left\|\re^{-tA_2}u_{0\varepsilon}+\int_0^t\re^{-(t-s)A_2}v_{\varepsilon}\rd s\right\|_{W^{1,p}(\Omega)}\notag\\
	&\leqslant C\|u_{0\varepsilon}\|_{W^{1,p}(\Omega)}+C\int_0^t\|v_\varepsilon\|_{W^{1,p}(\Omega)}\rd s
	\leqslant C\|A_1^\beta v_\varepsilon\|_p^{\delta_3}+C.
\end{align}
Plugging \eqref{v*} and \eqref{u*} into \eqref{vA1} and combining the fact that $v_{0\varepsilon}\in C_0^\infty(\Omega)\subset D(A_1)$ and $\|f_\varepsilon\|_{L^\infty(\Omega)}\leqslant C$, therefore we have
$$M_{1\varepsilon}(T)\leqslant CM_{1\varepsilon}^{\delta_3}(T)+C.$$
By utilizing the Young inequality, the proof of Lemma \ref{LemA1} is now finished.
\end{proof}

The boundedness of each $f_\varepsilon$ allows us, through heat semigroup estimates, to obtain a time-independent estimate for $\Theta_\varepsilon$ within a topological framework that is consistent with the one appearing in \eqref{Tmax}.
\begin{lemma}\label{LemA3}
Let $q>1$ and $A_3$ be the realization of $-\varepsilon\Delta+\lambda_0 I$ in $L^q(\Omega; \mathbb{R}^n)$ with domain $D(A_3)=\left\{\psi \in W^{2,q}(\Omega; \mathbb{R}^n) \mid \psi = \frac{\partial\psi}{\partial\nu}= 0 \text{ on } \partial\Omega\right\}$, where $\lambda_0$ is a positive constant.
Suppose $T_{\text{max}, \varepsilon} < \infty$ for some $\varepsilon \in (0,1)$. Then, for any $\beta \in (\frac12,1)$, there exists a constant $C(\sigma, \varepsilon) > 0$ such that the corresponding fractional power satisfies
$$\|A_3^\beta \Theta_\varepsilon(\cdot, t)\|_{L^q(\Omega)}\leqslant C(\sigma, \varepsilon),~\quad\forall~t \in (0, T_{\max,\varepsilon}).$$
\end{lemma}
\begin{proof}
Using a Duhamel representation associated with the third equation in \eqref{equ1*}, we get
\begin{align*}
	\Theta_{\varepsilon t}+ A_3\Theta_{\varepsilon} =h_\varepsilon:= \gamma_\varepsilon(\Theta_\varepsilon)|\nabla v_{\varepsilon}|^2 - f_\varepsilon(\Theta_\varepsilon)\nabla v_{\varepsilon}+\Theta_\varepsilon.
\end{align*}
By means of a standard smoothing estimate for the semigroup $(\re^{-tA})_{t\geqslant0}$, we obtain that
\begin{align}\label{Theta*}
	\|A_3^\beta\Theta_\varepsilon\|_q
	&=\left\|\re^{-tA_3}A_3^\beta\Theta_{0\varepsilon}+\int_0^tA_3^\beta\re^{-(t-s)A_3}h_{\varepsilon}\rd s\right\|_q\notag\\
	&\leqslant\|\re^{-tA_3}A_3^\beta\Theta_{0\varepsilon}\|_q+\int_0^tA_3^\beta\re^{-(t-s)A_3}\|h_{\varepsilon}\|_q\rd s\notag\\
	&\leqslant C\|A_3^\beta\Theta_{0\varepsilon}\|_q+C_\beta\int_0^t(t-s)^{-\beta}\|h_{\varepsilon}\|_q\rd s.
\end{align}
Let us set $$M_{3\varepsilon}(T):=\sup\limits_{t\in(0,T)}\|A_3^\beta \Theta_\varepsilon(\cdot,t)\|_q,~~\quad\text{for}~T\in(0,T_{\max,\varepsilon}).$$
The continuous embedding $W^{2,p}(\Omega)\hookrightarrow W^{1,\infty}(\Omega)$ which holds for some $p>N$, together with the boundedness of $\gamma_\varepsilon$, the Minkowski inequality, the H\"older inequality, the Gagliardo-Nirenberg inequality, the uniform bound $(f_\varepsilon)_{\varepsilon\in(0,1)}\subset L^\infty([0,\infty))$, \eqref{v*} and Lemma \ref{LemA1}, collectively ensure that
\begin{align*}
	\|h_\varepsilon\|_p&=\left\|\gamma_\varepsilon(\Theta_\varepsilon)|\nabla v_{\varepsilon}|^2 - f_\varepsilon(\Theta_\varepsilon)\nabla v_{\varepsilon}+\Theta_\varepsilon\right\|_q\\
	&\leqslant C\|v_\varepsilon\|_{W^{1,\infty}(\Omega)}^2+C\|f_\varepsilon\|_{L^\infty(\Omega)}\|v_\varepsilon\|_{W^{1,\infty}(\Omega)}+\|\Theta_\varepsilon\|_q\\
	&\leqslant C\|v_\varepsilon\|_{W^{2,p}(\Omega)}^2+C\|f_\varepsilon\|_{L^\infty(\Omega)}\|v_\varepsilon\|_{W^{2,p}(\Omega)}+C\|\Theta_\varepsilon\|_{W^{2\beta,q}(\Omega)}^\delta\|\Theta_\varepsilon\|_{L^1}^{1-\delta}+C\|\Theta_\varepsilon\|_{L^1}\\
	&\leqslant C\|v_\varepsilon\|_{W^{2,p}(\Omega)}^2+C\|\Theta_\varepsilon\|_{W^{2\beta,q}(\Omega)}^\delta\|\Theta_\varepsilon\|_{L^1}^{1-\delta}+C\|\Theta_\varepsilon\|_{L^1}\\
	&\leqslant C\|v_\varepsilon\|_{W^{2,p}(\Omega)}^2+C\|A_3^\beta\Theta_\varepsilon\|_q^\delta+C\\
	&\leqslant CM_{3\varepsilon}^{\delta}(T)+C,
\end{align*}
where $\delta=\frac{N(q-1)}{2\beta q+N(q-1)}<1$.
Plugging the above inequality into \eqref{Theta*}, with the help of the fact $\Theta_{0\varepsilon}\in C_0^\infty(\Omega)\subset D(A_3)$, we have
$$M_{3\varepsilon}(T)\leqslant CM_{3\varepsilon}^{\delta}(T)+C.$$
Thus we get $M_{3\varepsilon}(T)\leqslant C$,~for all~$T\in(0,T_{\max,\varepsilon})$ by using the Young inequality directly.
\end{proof}

The above lemmas eliminate the possibility of finite-time blow-up via the second alternative in \eqref{Tmax} for any such $\varepsilon\in(0,1)$. Consequently, we are now in a position to establish the global in time existence of solutions to the regularized problem, which is stated in the following lemma.
\begin{lemma}\label{Lemglobal}
For every $\varepsilon \in (0, 1)$, the solution to problem \eqref{equ1*} is global in time. 
That is, $T_{\text{max}, \varepsilon} = \infty$ in Lemma \ref{Lemlocal}.
\end{lemma}
\begin{proof}
Assume $T_{\max,\varepsilon}<\infty$ for each $\varepsilon\in(0,1)$, fixing $\eta\in(0,\frac14)$, we may pick $p\geqslant2$ sufficiently large and $\beta\in(\frac12,\frac34)$ such that $4\beta-\frac Np>2+2\eta$. 
By the property of the fractional powers established in Lemma \ref{LemA1}, this choice would ensure that the following embedding
$$\|v_\varepsilon\|_{W^{4\beta,p}(\Omega)}\hookrightarrow\|v_\varepsilon\|_{W^{2+2\eta,\infty}(\Omega)}.$$
Therefore one can find $C>0$ satisfies 
\begin{align}\label{v*eta}
	\|v_\varepsilon\|_{W^{2+2\eta,\infty}(\Omega)}\leqslant C.
\end{align}
In a similar way, we can take $q>1$ large enough and $\beta\in(\frac12,1)$ such that $2\beta-\frac Nq>1+\eta$, which implies the embedding  $$\|\Theta_\varepsilon\|_{W^{1+\eta,\infty}(\Omega)}\hookrightarrow\|v_\varepsilon\|_{W^{2\beta,q}(\Omega)}$$
holds. By the property of the fractional powers established in Lemma \ref{LemA3}, there exists a constant $C>0$ fulfilling
\begin{align}\label{Theta*eta}
	\|\Theta_\varepsilon\|_{W^{1+\eta,\infty}(\Omega)}\leqslant C.
\end{align}
By the second equation in \eqref{equ1*}, \eqref{v*eta} and the fact that $u_{0\varepsilon}\in C_0^\infty(\Omega)$, we get
\begin{align}\label{u*eta}
	\|u_\varepsilon\|_{W^{1+\eta,\infty}(\Omega)}
	&=\left\|\re^{-tA_2}u_{0\varepsilon}+\int_0^t\re^{-(t-s)A_2}v_{\varepsilon}\rd s\right\|_{W^{1+\eta,\infty}(\Omega)}\notag\\
	&\leqslant C\|u_{0\varepsilon}\|_{W^{1+\eta,\infty}(\Omega)}+C\int_0^t\|v_\varepsilon\|_{W^{1+\eta,\infty}(\Omega)}\rd s\notag\\
	&\leqslant C\|u_{0\varepsilon}\|_{W^{1+\eta,\infty}(\Omega)}+C\|v_\varepsilon\|_{W^{2+2\eta,\infty}(\Omega)}\leqslant C,
\end{align}
where $A_2=-\varepsilon\Delta(\cdot)$ in $L^p(\Omega; \mathbb{R}^n)$ with $D(A_2)=\left\{ \psi \in W^{2,p}(\Omega; \mathbb{R}^n) \mid \psi = 0 \text{ on } \partial\Omega\right\}$.	
If $T_{\max,\varepsilon}<\infty$, \eqref{v*eta}, \eqref{Theta*eta} and \eqref{u*eta} would contradict \eqref{Tmax}. Therefore we obtain $T_{\max,\varepsilon}=\infty$.
\end{proof}

\section{Further $\varepsilon$-independent estimates}
To construct the limit object $(v,u,\Theta)$ by suitably extracting subsequences of solutions to \eqref{equ1*}, we shall establish regularity properties that are uniform with respect to $\varepsilon$. As a first step in this direction, we present the following useful lemma, which relies on the Gagliardo-Nirenberg interpolation inequality. 
\begin{lemma}\label{Lemq}
Let $q\in\mathbb{R}$ and $p\in(1-\frac2N,2]\cap(0,2]$. Then for all $t>0$ and $\varepsilon\in(0,1)$, there exists $C=C(p,q)>0$ such that
	\begin{align}\label{GNqvare}
		\int_{\Omega}\left(1+\Theta_\varepsilon\right)^q\rd x
		\leqslant C\left(\int_{\Omega}(1+\Theta_\varepsilon)^{p-2} |\nabla\Theta_{\varepsilon}|^2\rd x\right)^{\frac{N(q-1)}{2+N(p-1)}}+C.
	\end{align}
Moreover, if $q<p+\frac 2N$, then for all $\varepsilon_1>0$, there exists $C=C(\varepsilon_1,p,q)>0$ such that the following inequality holds
\begin{align}\label{GNq2}
	\int_{\Omega}\left(1+\Theta_\varepsilon\right)^q\rd x
	\leqslant \varepsilon_1\int_{\Omega}(1+\Theta_\varepsilon)^{p-2} |\nabla\Theta_{\varepsilon}|^2\rd x+C.
\end{align}
\end{lemma}
\begin{proof}
If $q>1$, the Gagliardo-Nirenberg inequality ensures that
\begin{align*}
	\|\psi\|_{L^{\frac{2q}{p}}{\Omega}}^{\frac{2q}{p}}&\leqslant C\|\psi\|_{W^{1,2}(\Omega)}^{\frac{2q}{p}\delta}\cdot\|\psi\|_{L^{\frac2p}}^{\frac{2q}{p}(1-\delta)}
	\leqslant C\|\nabla\psi\|_{L^{2}(\Omega)}^{\frac{2q}{p}\delta}\cdot\|\psi\|_{L^{\frac2p}}^{\frac{2q}{p}(1-\delta)}+C\|\psi\|_{L^{\frac2p}}^{\frac{2q}{p}},
\end{align*}
where $\frac{2q}{p}\delta=\frac{2N(q-1)}{2+N(p-1)}>0$. Let $\psi=(1+\Theta_\varepsilon)^{\frac p2}$ in the above inequality, together with Lemma \ref{Lem-energy*}. Then we get
\begin{align*}
		\int_{\Omega}\left(1+\Theta_\varepsilon\right)^q\rd x
		&\leqslant C\left(\int_{\Omega}(1+\Theta_\varepsilon)^{p-2} |\nabla\Theta_{\varepsilon}|^2\rd x\right)^{\frac{q}{p}\delta}\|1+\Theta_\varepsilon\|_{L^1(\Omega)}^{q(1-\delta)}+C\|1+\Theta_\varepsilon\|_{L^1(\Omega)}^{q}\\
		&\leqslant C\left(\int_{\Omega}(1+\Theta_\varepsilon)^{p-2} |\nabla\Theta_{\varepsilon}|^2\rd x\right)^{\frac{N(q-1)}{2+N(p-1)}}+C.
\end{align*}
Further, if $q<p+\frac 2N$, we can obtain $\frac{N(q-1)}{2+N(p-1)}<1$. Then, from the above inequality, Lemma \ref{Lem-energy*} and the Young inequality, it yields
\begin{align*}
	\int_\Omega(1+\Theta_\varepsilon)^q\rd x
	\leqslant\varepsilon_1 \int_{\Omega}(1+\Theta_\varepsilon)^{p-2} |\nabla\Theta_{\varepsilon}|^2\rd x+C_{\varepsilon_1}.
\end{align*}
If $q\leqslant1$, by using Lemma \ref{Lem-energy*} and the Young inequality, the following inequality
$$\int_{\Omega}(1+\Theta_\varepsilon)^q\rd x\leqslant C,$$
holds directly. Thus, the proof of Lemma \ref{Lemq} is now finished.
\end{proof}

In order to extract strongly convergent subsequences of $\Theta_\varepsilon$, a crucial step is to obtain uniform estimates for its spatial gradient. The following lemma, which provides an $\varepsilon$-independent bound for $\nabla\Theta_\varepsilon$ in $L^r$ with $r\in\big[1, \frac{N+2}{N+1}\big)$, will be essential in handling the nonlinear coupling terms in the limit procedure.
\begin{lemma}\label{LemThetanabla}
Let $r\in\big[1,\frac{N+2}{N+1}\big)$ and $T>0$. Then for all $\varepsilon\in(0,1)$, there exists $C(r,T)>0$ such that
\begin{equation*}
 	\int_{0}^{T} \int_{\Omega} |\nabla\Theta_{\varepsilon}|^r\rd x\rd t 
 	\leqslant C(r,T).
\end{equation*}
\end{lemma}
\begin{proof}
Employing the third equation in \eqref{equ1*} together with the homogeneous Neumann boundary conditions for $\Theta_\varepsilon$, via an integration by parts, we obtain
\begin{align*}
	\frac1p\frac{\rd}{\rd t}\int_{\Omega}(1+\Theta_\varepsilon)^{p}\rd x
	&=\int_{\Omega}(1+\Theta_\varepsilon)^{p-1}\left(\Delta\Theta_{\varepsilon}+\gamma_\varepsilon(\Theta_\varepsilon)|\nabla v_{\varepsilon}|^2-f_\varepsilon(\Theta_\varepsilon)\nabla v_{\varepsilon}\right)\rd x\notag\\
	&=(1-p)\int_{\Omega}(1+\Theta_\varepsilon)^{p-2}\cdot |\nabla\Theta_{\varepsilon}|^2\rd x+\int_{\Omega}\gamma_\varepsilon(\Theta_\varepsilon)|\nabla v_{\varepsilon}|^2\cdot(1+\Theta_\varepsilon)^{p-1}\rd x\notag\\
	&~~~-\int_{\Omega}f_\varepsilon(\Theta_\varepsilon)\nabla v_{\varepsilon}\cdot(1+\Theta_\varepsilon)^{p-1}\rd x.
\end{align*}
Let $p\in(0,1)\cap(1-\frac2N,1)$ in the above equality. By using the Young inequality, \eqref{gamma*bound}, \eqref{falpha} and \eqref{GNqvare}, then we have
\begin{align}\label{ThetaYoung}
	-\frac1p\frac{\rd}{\rd t}\int_{\Omega}(1&+\Theta_\varepsilon)^{p}\rd x
	+(1-p)\int_{\Omega}(1+\Theta_\varepsilon)^{p-2}\cdot|\nabla \Theta_{\varepsilon}|^2\rd x\notag\\
	&=-\int_{\Omega}\gamma_\varepsilon(\Theta_\varepsilon)|\nabla v_{\varepsilon}|^2\cdot(1+\Theta_\varepsilon)^{p-1}\rd x
	+\int_{\Omega}f_\varepsilon(\Theta_\varepsilon)\nabla v_{\varepsilon}\cdot(1+\Theta_\varepsilon)^{p-1}\rd x\notag\\
	&\leqslant-k_\gamma\int_{\Omega}|\nabla v_{\varepsilon}|^2\cdot(1+\Theta_\varepsilon)^{p-1}\rd x
	+\int_{\Omega}\left(C_{\gamma}K_f^2|f_\varepsilon(\Theta_\varepsilon)|^2+k_\gamma |\nabla v_{\varepsilon}|^2\right)\cdot(1+\Theta_\varepsilon)^{p-1}\rd x\notag\\
	&\leqslant C\int_{\Omega}(1+\Theta_\varepsilon)^{p+2\alpha-1}\rd x
	\leqslant \varepsilon_1\int_{\Omega}(1+\Theta_\varepsilon)^{p-2}|\nabla\Theta_{\varepsilon}|^2\rd x+C_{\varepsilon_1}.
\end{align}
Here we have used $p+2\alpha-1<p+\frac 2N$. Taking $\varepsilon_1$ small enough such that $1-p-\varepsilon_1>0$, from \eqref{ThetaYoung}, the Young inequality with $p<1$ and Lemma \ref{Lem-energy*}, we obtain

\begin{align}\label{Theta0T}
	\int_0^T\int_{\Omega}(1+\Theta_\varepsilon)^{p-2}\cdot|\nabla \Theta_{\varepsilon}|^2\rd x\rd s
	&\leqslant C\left(\int_{\Omega}(1+\Theta_\varepsilon(\cdot,T))^{p}\rd x-\int_{\Omega}(1+\Theta_{0\varepsilon})^p\right)+C\notag\\
	&\leqslant C\int_{\Omega}\Theta_{\varepsilon}(\cdot,T)\rd x+C\leqslant C(p,T),
\end{align}
for all $T>0$ and $\varepsilon\in(0,1)$.
Combining the above inequality and using the Young inequality again, it yields that
\begin{align*}
	\int_{0}^{T} \int_{\Omega} |\nabla\Theta_{\varepsilon}|^r\rd x\rd t 
	&=	\int_{0}^{T} \int_{\Omega} \left((1+\Theta_\varepsilon)^{p-2}|\nabla\Theta_{\varepsilon}|^2\right)^{\frac r2}\cdot(1+\Theta_\varepsilon)^{\frac{r(2-p)}{2}}\rd x\rd t \\
	&\leqslant\int_0^T\int_{\Omega}(1+\Theta_\varepsilon)^{p-2}\cdot |\nabla\Theta_{\varepsilon}|^2\rd x\rd t
	+\int_0^T\int_{\Omega}(1+\Theta_\varepsilon)^{\frac{r(2-p)}{2-r}}\rd x\rd t,
\end{align*}
where $\frac{r(2-p)}{2-r}<p+\frac 2N$. From \eqref{Theta0T} and Lemma \ref{Lemq}, we deduce that
\begin{align*}
	\int_{0}^{T} \int_{\Omega} |\nabla\Theta_{\varepsilon}|^r\rd x\rd t 
	\leqslant C\int_0^T\int_{\Omega}(1+\Theta_\varepsilon)^{p-2}\cdot |\nabla\Theta_{\varepsilon}|^2\rd x\rd t+C
	\leqslant C(r,T),
\end{align*}
for all $T>0$ and $\varepsilon\in(0,1)$. Thus, the claim is shown.
\end{proof}

The following inequality is established directly from Lemma \ref{Lemq} and \eqref{Theta0T}.
\begin{remark}\label{RemarkThetaq}
Let $T>0$ and $q\in\big[1,\frac{N+2}{N}\big)$. Then, for all $\varepsilon\in(0,1)$, there exists $C=C(q,T)>0$ such that  
\begin{align*}
	\int_0^T\int_{\Omega}\left(1+\Theta_\varepsilon\right)^q\rd x\rd s\leqslant C.
\end{align*}
\end{remark}

To pass to the limit as $\varepsilon \searrow0$, we need compactness properties for the sequences ${v_\varepsilon}$, ${u_\varepsilon}$ and ${\Theta_\varepsilon}$. The following lemma, which collects estimates on their time derivatives in various dual spaces, will allow us to apply the Aubin–Lions compactness theorem and thereby obtain strong convergence of the relevant quantities.
\begin{lemma}\label{Lem*vuTheta}
Let $\lambda>N+2$ and $T>0$. For all $\varepsilon\in(0,1)$,  then there exist $C(T)>0$ and $C(\lambda,T)>0$ such that
	$$\int_{0}^{T} \int_{\Omega}|\nabla v_{\varepsilon}|^2\rd x\rd t+\int_{0}^{T} \left\| v_{\varepsilon t}(\cdot,t) \right\|_{\big(W^{2,2}_{0}(\Omega)\big)^{*}}^2 \rd t+\int_0^T\int_\Omega u_{\varepsilon t}^2\rd x\rd t\leqslant C(T),$$
and that
$$\int_{0}^{T} \left\| \Theta_{\varepsilon t}(\cdot,t) \right\|_{\big(W^{1,\lambda}(\Omega)\big)^{*}} \rd t \leq C(\lambda,T).$$
\end{lemma}
\begin{proof}
Integrating \eqref{energy*v} over $(0,T)$, dropping some negative terms and using the Cauchy-Schwarz inequality, together with \eqref{gamma*bound} and \eqref{falpha}, 
it yields
\begin{align*}
	k_\gamma\int_0^T\int_{\Omega}|\nabla v_{\varepsilon}|^2\rd x\rd t
	&\leqslant\int_0^T\int_{\Omega}f_\varepsilon(\Theta_\varepsilon)\nabla v_{\varepsilon}\rd x\rd t\\
	&\leqslant\frac{k_\gamma}{2}\int_0^T\int_{\Omega}|\nabla v_{\varepsilon}|^2\rd x\rd t
	+C_{k_\gamma}K_f^2\int_0^T\int_{\Omega}(1+\Theta_\varepsilon)^{2\alpha}\rd x\rd t.
\end{align*}
From Remark \ref{RemarkThetaq}, there exists $C(T)>0$ such that
\begin{align}\label{nablav*}
	\int_0^T\int_{\Omega}|\nabla v_{\varepsilon}|^2\rd x\rd t\leqslant C(T).
\end{align}
For fixed $\psi\in W_0^{2,2}(\Omega)$ satisfying $\|\psi\|_2+\|\nabla\psi\|_2+\|\Delta\psi\|_2\leqslant1$,  multiplying the first equation in \eqref{equ1*} by $\psi$, integrating over $\Omega\times(0,T) (T>0)$ and using the H\"older inequality, we obtain
\begin{align*}
	\left|\int_\Omega v_{\varepsilon t}\psi\rd x\right|
	&\leqslant\int_\Omega \varepsilon|\Delta v_{\varepsilon}|\cdot|\Delta\psi|\rd x+\int_\Omega \gamma_\varepsilon(\Theta_\varepsilon)|\nabla v_{\varepsilon}|\cdot|\nabla\psi|\rd x\notag\\
	&~~~+a\int_\Omega |\nabla u_{\varepsilon}|\cdot|\nabla\psi|\rd x+\int_\Omega |f_\varepsilon(\Theta_\varepsilon)|\cdot|\nabla\psi|\rd x\\
	&\leqslant\varepsilon\|\Delta v_\varepsilon\|_2+K_\gamma\|\nabla v_\varepsilon\|_2
	+a\|\nabla u_\varepsilon\|_2+\|f_\varepsilon(\Theta_\varepsilon)\|_2.
\end{align*}
From \eqref{falpha}, \eqref{nablav*}, Lemma \ref{Lem-energy*} and Remark \ref{RemarkThetaq}, together with the above inequality, we have
\begin{align*}
	\int_{0}^{T} \left\| v_{\varepsilon t}(\cdot,t) \right\|_{\big(W^{2,2}_{0}(\Omega)\big)^{*}}^2 \rd t
	&\leqslant C\varepsilon^2\int_0^T\int_\Omega|\Delta v_\varepsilon|^2\rd x\rd t
	+C\int_0^T\int_\Omega|\nabla v_\varepsilon|^2\rd x\rd t\\
	&~~~+C\int_0^T\int_\Omega|\nabla u_\varepsilon|^2\rd x\rd t
	+CK_f^2\int_0^T\int_\Omega(1+\Theta_\varepsilon)^{2\alpha}\rd x\rd t
	\leqslant C(T),
\end{align*}
where $C(T)>0$ is a constant depending on $T$, but independent of the parameter $\varepsilon$.
By using the second equation in \eqref{equ1*}, the Cauchy-Schwarz inequality and Lemma \ref{Lem-energy*}, it yields
\begin{align*}
	\int_0^T\int_\Omega u_{\varepsilon t}^2\rd x\rd t
	&=\int_0^T\int_\Omega\left(\varepsilon\Delta u_\varepsilon+v_{\varepsilon}\right)^2\rd x\rd t\\
	&\leqslant2\varepsilon^2\int_0^T\int_\Omega|\Delta u_\varepsilon|^2\rd x\rd t
	+2\int_0^T\int_\Omega v_\varepsilon^2\rd x\rd t
	\leqslant C(T).
\end{align*}
Multiplying the third equation in \eqref{equ1*} by $\psi$, integrating over $\Omega\times(0,T) (T>0)$ by parts, and using the H\"older inequality, \eqref{gamma*bound} and Sobolev embedding $W^{1,\lambda}(\Omega)\hookrightarrow L^\infty(\Omega)$, for fixed $\psi\in W^{1,\lambda}(\Omega)$ fulfilling $\|\psi\|_{W^{1,\lambda}(\Omega)}\leqslant1,~\lambda>N+2$, we obtain 
\begin{align*}
	\left|\int_\Omega\Theta_{\varepsilon t}\psi\rd x\right|
	&=\left|-\int_\Omega\nabla\Theta_\varepsilon\nabla\psi\rd x+\int_\Omega\gamma_\varepsilon(\Theta_\varepsilon)|\nabla v_{\varepsilon}|^2\psi\rd x-\int_\Omega f_\varepsilon(\Theta_\varepsilon)\nabla v_{\varepsilon}\psi\rd x\right|\\
	&\leqslant\|\nabla\Theta_\varepsilon\|_{L^{\frac{\lambda}{\lambda-1}}}\|\nabla\psi\|_{L^\lambda(\Omega)}+K_\gamma\|\psi\|_{L^\infty(\Omega)}\|\nabla v_\varepsilon\|_2^2+\|\psi\|_{L^\infty(\Omega)}\|f_\varepsilon(\Theta_\varepsilon)\|_2\|\nabla v_\varepsilon\|_2\\
	&\leqslant\|\nabla\Theta_\varepsilon\|_{L^{\frac{\lambda}{\lambda-1}}}\|\nabla\psi\|_{L^\lambda(\Omega)}+K_\gamma\|\psi\|_{W^{1,\lambda}(\Omega)}\|\nabla v_\varepsilon\|_2^2+\|\psi\|_{W^{1,\lambda}(\Omega)}\|f_\varepsilon(\Theta_\varepsilon)\|_2\|\nabla v_\varepsilon\|_2,
\end{align*}
where $\frac{\lambda}{\lambda-1}<\frac{N+2}{N+1}$ since the assumption $\lambda>N+2$. Recalling the fact that $\|\psi\|_{W^{1,\lambda}(\Omega)}\leqslant1$, \eqref{falpha}, Lemma \ref{Lemq}, Lemma \ref{LemThetanabla}, Remark \ref{RemarkThetaq} and \eqref{nablav*}, together with the Cauchy inequality, it yields that
\begin{align*}
	\int_0^T\|\Theta_{\varepsilon t}(\cdot,t)\|_{\big(W^{1,\lambda}(\Omega)\big)^*}\rd t
	\leqslant C(\lambda,T),
\end{align*}
for all $T>0$ and $\varepsilon\in(0,1)$.
Thus, the proof of Lemma \ref{Lem*vuTheta} is now finished.
\end{proof}

With the $\varepsilon$-independent estimates established in the previous lemmas, we are now in a position to extract convergent subsequences and identify their limits. The following lemma summarizes the convergence properties obtained through a standard diagonal extraction procedure, which will serve as the foundation for passing to the limit in the regularized system.
\begin{lemma}\label{Lemwekav}
There exists $\{\varepsilon_k\}_{k\in\mathbb{N}}\subset(0,1)$  such that $\varepsilon=\varepsilon_k\searrow0,~k\rightarrow\infty$ and that with some functions 
\begin{equation*}
	\begin{cases}
		v\in L^\infty\left((0,\infty); L^2(\Omega)\right)\cap	L_{loc}^2\left([0,\infty); W_0^{1,2}(\Omega)\right)\cap L_{loc}^2\left([0,\infty)\times\Omega\right),\\
		u\in C\left([0,T],L^2(\Omega)\right)\cap L^\infty\left([0,\infty); W_0^{1,2}(\Omega)\right),\\
		\Theta\in\bigcap_{q\in[1,\frac{N+2}{N})} L_{loc}^q\left([0,\infty)\times\Omega\right)\cap \bigcap_{r\in(1,\frac{N+2}{N+1})}L_{loc}^r\left([0,\infty); W^{1,r}(\Omega)\right),
	\end{cases}
\end{equation*}
and $$f(\Theta)\in L_{loc}^2\left([0,\infty)\times\Omega\right),\quad \gamma(\Theta)\in L_{loc}^2\left([0,\infty)\times\Omega\right),$$
satisfying $\Theta\geqslant0$ a.e. in $\Omega\times(0,\infty)$, we have
\begin{align}\label{weak*vuTheta}
	v_\varepsilon&\rightarrow v\quad\quad\text{a.e.~in}~\Omega\times(0,\infty)~\text{and~in}~L_{loc}^2\left([0,\infty)\times\Omega\right),\notag\\
	\nabla v_{\varepsilon}&\rightharpoonup \nabla v\quad~\text{in}~L_{loc}^2\left([0,\infty)\times\Omega\right),\notag\\
	u_\varepsilon&\rightarrow u~~~\quad\text{in}~C\left([0,T],L^2(\Omega)\right),\notag\\
	\nabla u_{\varepsilon}&\rightharpoonup\nabla u\quad~\text{in}~L_{loc}^2\left([0,\infty)\times\Omega\right),\notag\\
	\Theta_\varepsilon&\rightarrow \Theta\quad\quad\text{a.e.~in}~\Omega\times(0,\infty)~\text{and~in}~ L_{loc}^q\left([0,\infty)\times\Omega\right),\\
	\nabla\Theta_\varepsilon&\rightharpoonup\nabla\Theta~~~~ \text{in}~L_{loc}^r\left([0,\infty)\times\Omega
	\right),\notag\\
	f_\varepsilon(\Theta_\varepsilon)&\rightarrow f(\Theta)\quad\text{in}~L_{loc}^2\left([0,\infty)\times\Omega\right),\notag\\
	\gamma_\varepsilon(\Theta_\varepsilon)&\rightarrow
	\gamma(\Theta)\quad\text{in}~L_{loc}^2\left([0,\infty)\times\Omega\right)\notag,
\end{align}
where $q\in\big[1,\frac{N+2}{N}\big)$ and $r\in\left(1,\frac{N+2}{N+1}\right)$.
Moreover, we have
\begin{align}\label{limitu}
	-\int_{0}^{\infty} \int_{\Omega} &v \varphi_t\rd x\rd t
	- \int_{\Omega} u_{0t} \varphi(\cdot, 0)\rd x\notag\\
	&=-\int_{0}^{\infty} \int_{\Omega} \gamma(\Theta)\nabla v \nabla\varphi \rd x\rd t
	-a \int_{0}^{\infty} \int_{\Omega} \nabla u \nabla\varphi \rd x\rd t
	+ \int_{0}^{\infty} \int_{\Omega} f(\Theta) \nabla\varphi\rd x\rd t,
\end{align}
for all $\varphi\in C_0^\infty\left([0,\infty)\times\Omega\right)$, and 
\begin{equation}\label{vuta.e.}
	u_t=v,~\quad\text{a.e.~in}~\Omega\times(0,\infty).
\end{equation}
\end{lemma}
\begin{proof}
\textbf{Step 1.
We shall give firstly the proof of \eqref{weak*vuTheta} in sequence.} 
According to Lemma \ref{Lem*vuTheta}, for all $\varepsilon\in(0,1)$ and $T>0$, we obtain
$$v_\varepsilon\in L^2\left(0,T; W_0^{1,2}(\Omega)\right)\quad \text{and}\quad v_{\varepsilon t}\in L^2\left(0,T; \big(W_0^{2,2}(\Omega)\big)^*\right).$$
By using Sobolev embedding $W_0^{1,2}(\Omega)\hookrightarrow\hookrightarrow L^2(\Omega)\hookrightarrow  \big(W_0^{2,2}(\Omega)\big)^*$, together with the Aubin-Lions Lemma, we get
\begin{equation}\label{vL2}
	v_\varepsilon\rightarrow v\quad\text{in}~L_{loc}^2\left([0,\infty)\times\Omega\right),
\end{equation}
and
$$\nabla v_\varepsilon\rightharpoonup \nabla v~~\text{in}~L_{loc}^2\left([0,\infty)\times\Omega\right).$$
Then $v_\varepsilon\rightarrow v$ a.e. in $(0,\infty)\times\Omega$ is obtained from \eqref{vL2}.
Combining Lemma \ref{Lem-energy*} and Lemma \ref{Lem*vuTheta}, for all $\varepsilon\in(0,1)$ and $T>0$, we get
$$u_\varepsilon\in L^\infty\left(0,T; W_0^{1,2}(\Omega)\right)\quad\text{and}\quad u_{\varepsilon t}\in L^2\left(0,T; L^2(\Omega)\right).$$
Using the Aubin-Lions Lemma again with Sobolev embedding $W_0^{1,2}(\Omega)\hookrightarrow\hookrightarrow L^2(\Omega)\hookrightarrow L^2(\Omega)$ and the weak compactness for $L^2(\Omega)$, it yields that
$$u_\varepsilon\rightarrow u\quad\text{in}~C\left([0,T], L^2(\Omega)\right)\quad\text{and}\quad\nabla u_{\varepsilon}\rightharpoonup \nabla u\quad\text{in}~L_{loc}^2\left([0,\infty)\times\Omega\right).$$
By Lemma \ref{LemThetanabla} and Lemma \ref{Lem*vuTheta} again, for all $\varepsilon\in(0,1)$ and $T>0$, we have 
$$\Theta_\varepsilon\in L^r\left(0,T; W^{1,r}(\Omega)\right)\quad\text{and}\quad \Theta_{\varepsilon t}\in L^1\left(0,T; \big(W^{1,\lambda}(\Omega)\big)^*\right),$$
where $r\in\big[1,\frac{N+2}{N+1}\big)$. By the weak compactness for $L^r(\Omega) \big(1<r<\frac{N+2}{N+1}\big)$ and the Aubin-Lions Lemma again with Sobolev embedding $W^{1,r}(\Omega)\hookrightarrow\hookrightarrow L^1(\Omega)\hookrightarrow (W^{1,\lambda}(\Omega))^*$, we get
$$\nabla\Theta_\varepsilon\rightharpoonup\nabla\Theta~~\text{in}~L_{loc}^r\left([0,\infty)\times\Omega\right)\quad\text{and}\quad\Theta_\varepsilon\rightarrow\Theta~~\text{in}~L^r\left(0,T,L^1(\Omega)\right),$$
where $r\in\big(1,\frac{N+2}{N+1}\big)$. Then we have $\Theta_\varepsilon\rightarrow\Theta$ in $L^1\left([0,T)\times\Omega\right)$ for all $T>0$. Consequently, we obtain
\begin{equation}\label{Theta.e.}
	\Theta_\varepsilon\rightarrow\Theta\quad\text{a.e.~in}~(0,\infty)\times\Omega.
\end{equation}
Combining \eqref{Theta.e.} and Remark \ref{RemarkThetaq}, together with the Vitali Convergence Theorem, we have directly
$$\Theta_\varepsilon\rightarrow\Theta\quad\text{in}~L_{loc}^q\left([0,\infty)\times\Omega\right),$$  
for all $q\in\left[1,\frac{N+2}{N}\right)$.
Since $\{\gamma_\varepsilon\}_{\varepsilon\in(0,1)}\subset C^\infty([0,\infty)),~\{f_\varepsilon\}_{\varepsilon\in(0,1)}\subset C^\infty([0,\infty))$, \eqref{gamma*bound}, \eqref{falpha}, \eqref{Theta.e.} and Lemma \ref{Lemq}, together with the Vitali Convergence Theorem, it is easy to see
$$f_\varepsilon(\Theta_\varepsilon)\rightarrow f(\Theta)\quad\text{in}~L_{loc}^2\left([0,\infty)\times\Omega\right),$$
and
$$ \gamma_\varepsilon(\Theta_\varepsilon)\rightarrow
\gamma(\Theta)\quad\text{in}~L_{loc}^2\left([0,\infty)\times\Omega\right).$$
Consequently, \eqref{weak*vuTheta} is obtained.

\textbf{Step 2.
Next, we will prove that \eqref{limitu} and \eqref{vuta.e.} hold through the limit process of $\varepsilon=\varepsilon_k\searrow0,~k\rightarrow\infty$.}
Multiplying the first equation in \eqref{equ1*} by $\varphi$, integrating over $\Omega\times(0,T) (T>0)$ by parts, for all $\varphi\in C_0^\infty\left([0,\infty)\times\Omega\right)$, we obtain
\begin{align*}
	-\int_{0}^{\infty} \int_{\Omega} v_\varepsilon \varphi_t\rd x\rd t
	&- \int_{\Omega} v_{0\varepsilon} \varphi(\cdot, 0)\rd x\\
	&=-\varepsilon\int_0^\infty\int_\Omega v_\varepsilon\Delta^2\varphi\rd x\rd t-\int_{0}^{\infty} \int_{\Omega} \gamma_\varepsilon(\Theta_\varepsilon)\nabla v_\varepsilon \nabla\varphi \rd x\rd t\\
	&~~~- a \int_{0}^{\infty} \int_{\Omega} \nabla u_\varepsilon \nabla\varphi \rd x\rd t
	+ \int_{0}^{\infty} \int_{\Omega} f_\varepsilon(\Theta_\varepsilon) \nabla\varphi\rd x\rd t.
\end{align*}
Together with \eqref{weak*vuTheta} and \eqref{smoothdata}, we get
\begin{align*}
	-\int_{0}^{\infty} \int_{\Omega} &v \varphi_t\rd x\rd t
	- \int_{\Omega} u_{0t} \varphi(\cdot, 0)\rd x\notag\\
	&=-\int_{0}^{\infty} \int_{\Omega} \gamma(\Theta)\nabla v \nabla\varphi \rd x\rd t
	-a \int_{0}^{\infty} \int_{\Omega} \nabla u \nabla\varphi \rd x\rd t
	+ \int_{0}^{\infty} \int_{\Omega} f(\Theta) \nabla\varphi\rd x\rd t.
\end{align*}
as $\varepsilon=\varepsilon_k\searrow0,~k\rightarrow\infty$. Then \eqref{limitu} holds.
Finally, multiplying the second equation in \eqref{equ1*} by $\varphi$, integrating over $\Omega\times(0,T) (T>0)$, we obtain
\begin{align*}
	-\int_0^\infty\int_\Omega u_\varepsilon\varphi_t\rd x\rd t=\varepsilon\int_0^\infty\int_\Omega u_\varepsilon\Delta\varphi\rd x\rd t+\int_0^\infty\int_\Omega v_\varepsilon\varphi\rd x\rd t,
\end{align*}
for all $\varphi\in C_0^\infty\left([0,\infty)\times\Omega\right)$. From \eqref{weak*vuTheta} and letting $\varepsilon=\varepsilon_k\searrow0,~k\rightarrow\infty$, \eqref{vuta.e.} holds.
Thus the proof of Lemma \ref{Lemwekav} is finished now.
\end{proof}

\section{Solution properties of $\Theta$}
To handle the time derivatives in the weak formulation, we shall employ the Steklov averaging technique \cite{Steklov1993}. The following lemma establishes a useful relation between the Steklov averages of $\widehat{v}$ and $\widehat{u}$, which will allow us to pass to the limit in terms involving $v$ and $\nabla v$.
\begin{lemma}\label{LemSteklov}
Let $h>0$ and $(v,u,\Theta)$ be as in Lemma \ref{Lemwekav}. We denote that
$$(S_h\varphi)(x,t):=\frac1h\int_{t-h}^h\varphi(x,s)\rd s,\quad x\in\Omega,~t\in\mathbb{R},~\varphi\in L_{loc}^1(\Omega\times\mathbb{R}).$$
Let
\begin{align}\label{hatv}
\widehat{v}(x,t) := 
\begin{cases} 
	v(x,t), & x \in \Omega,\ t > 0, \\
	u_{0t}(x), & x \in \Omega,\ t < 0,
\end{cases}
\end{align}
and 
\begin{equation}\label{hatu}
  \widehat{u}(x,t) := 
  \begin{cases} 
  	u(x,t), & x \in \Omega,\ t > 0, \\
  	u_{0}(x) + t u_{0t}(x), & x \in \Omega,\ t < 0.
  \end{cases}
\end{equation}
Then we have
\begin{align}\label{Suth}
	(S_h \nabla\widehat{v})(x,t) &= \frac{\nabla\widehat{u}(x,t) - \nabla\widehat{u}(x,t-h)}{h},
\end{align}
for a.e. $(x,t)\in\Omega\times\mathbb{R}~\text{and~each}~h>0$.
Moreover, we have
\begin{equation}\label{hatsteklov}
	S_h \widehat{v}\to v\quad\text{and}\quad S_h \nabla\widehat{v}\to \nabla v\quad \text{both~in}~L^2_{\text{loc}}(\Omega \times [0,\infty)),
\end{equation}
as $h=h_k\searrow0,~k\rightarrow\infty$.
\end{lemma}
\begin{proof}
A standard theory of Steklov averages [Lemma 3.2,~\ref{Steklov1993}], applied to the inclusions $\widehat{v} \in L^2_{loc}(\Omega \times \mathbb{R})$ and $\nabla\widehat{v}\in L^2_{loc}(\Omega \times \mathbb{R})$ from Lemma \ref{Lemwekav}, implies $S_h \widehat{v} \to \widehat{v}$ and $S_h \nabla\widehat{v} \to \nabla\widehat{v}$ both in $L^2_{loc}(\Omega \times \mathbb{R})$ as $h=h_k \searrow 0, k\rightarrow\infty$. 
In view of \eqref{hatv}, then \eqref{hatsteklov} follows.
Notice that
\begin{align*}
  \nabla\widehat{u}_t=\nabla\widehat{v},\quad\text{a.e.~in}~\Omega\times\mathbb{R},
\end{align*}
holds through the standard integration technique with \eqref{hatv} and \eqref{hatu}.
Next, we will prove \eqref{Suth}.
Recalling the definition of $(S_h \nabla\widehat{v})(x,t)$ and the fact that $$\{(s,t): t-h\leqslant s\leqslant t, -\infty\leqslant t\leqslant\infty\}=\{(t,s): s\leqslant t\leqslant s+h, -\infty\leqslant s\leqslant\infty\},$$
together with the Fubini Theorem, we obtain
\begin{align*}
   \int_{-\infty}^{\infty} \int_{\Omega} \varphi(x,t) \cdot(S_h \nabla\widehat{v})(x,t)\rd x\rd t 
   &=\frac1h\int_{-\infty}^{\infty} \int_{\Omega} \varphi(x,t) \cdot \left\{ \int_{t-h}^{t} \nabla\widehat{u}_t(x,s) \rd s \right\} \rd x\rd t \\
   &=\frac1h\int_{-\infty}^{\infty} \int_{\Omega}\left\{ \int_{s}^{s+h}\varphi(x,t) \rd t \right\}\cdot\nabla\widehat{u}_t(x,s) \rd x\rd s\\
   &=-\frac1h\int_{-\infty}^{\infty} \int_{\Omega}\frac{\rd }{\rd t}\left\{ \int_{s}^{s+h}\varphi(x,t) \rd t \right\}\cdot\nabla\widehat{u}(x,s) \rd x\rd s\\
   &=-\frac1h\int_{-\infty}^{\infty} \int_{\Omega}\varphi(x,s+h)\cdot\nabla\widehat{u}(x,s) \rd x\rd s\\
   &~~~+\frac1h\int_{-\infty}^{\infty} \int_{\Omega}\varphi(x,s)\cdot\nabla\widehat{u}(x,s) \rd x\rd s\\
   &=-\frac1h\int_{-\infty}^{\infty} \int_{\Omega}\varphi(x,s)\cdot\nabla\widehat{u}(x,s-h) \rd x\rd s\\
   &~~~+\frac1h\int_{-\infty}^{\infty} \int_{\Omega}\varphi(x,s)\cdot\nabla\widehat{u}(x,s) \rd x\rd s\\
   &=\int_{-\infty}^{\infty} \int_{\Omega}\varphi(x,s)\cdot\frac{\nabla\widehat{u}(x,s)-\nabla\widehat{u}(x,s-h)}{h} \rd x\rd s,
\end{align*}
for all $\varphi\in C_0^\infty(\Omega\times\mathbb{R})$, which implies \eqref{Suth} is valid.
\end{proof}

The next lemma provides the strong convergence of $\sqrt{\gamma_\varepsilon(\Theta_\varepsilon)}\nabla v_\varepsilon$ in $L^2\left(\Omega\times(0,T)\right)$, which is essential for passing to the limit in the quadratic nonlinearity $\gamma(\Theta)|\nabla v|^2$ appearing in the heat equation of \eqref{equ1*}.
\begin{lemma}\label{Lemgammato}
For any $T>0$, let $(v,u,\Theta)$ and $\{\varepsilon_{k}\}_{k\in\mathbb{N}}$ be as in Lemma \ref{Lemwekav}. Then we have
$$\sqrt{\gamma_\varepsilon(\Theta_\varepsilon) }\nabla v_{\varepsilon}
\to \sqrt{\gamma(\Theta) }\nabla v\quad\text{in}~L^2\left(\Omega\times(0,T)\right),$$
as $\varepsilon=\varepsilon_{k}\searrow0,~k\rightarrow\infty$.
\end{lemma}
\begin{proof}
\textbf{Step 1.}
By a standard approximation argument, we thus infer that \eqref{limitu} remains valid for all $\varphi\in L^2\left((0,\infty); W_0^{1,2}(\Omega)\right)$ such that $\varphi_t \in L^2((0, \infty)\times\Omega)$ and $\varphi = 0$ a.e. in $\Omega \times (T, \infty)$ for any $T > 0$. Consequently, for arbitrary $h > 0$, we may choose the test function $$\varphi(x, t) := \xi(t) \cdot (S_h \widehat{v})(x, t),\quad(x, t) \in \Omega \times (0, \infty),$$
where $\xi(t)\in C_0^\infty([0,\infty))$ such that $\xi(0)\equiv1$ on $[0,T) (T>0)$ and $\xi'(t)\leqslant0$ for all $t\geqslant0$.
Then, from Lemma \ref{LemSteklov}, it yields that
$$\displaystyle\varphi_t=\xi'(t)\cdot(S_h\widehat{v})(x,t)+\xi(t)\cdot\frac{\widehat v(x,t)-\widehat v(x,t-h)}{h},$$
and 
$$\displaystyle\nabla\varphi=\xi(t)\cdot(S_h\nabla\widehat{v})(x,t)=\xi(t)\cdot\frac{\nabla\widehat{u}(x,t)-\nabla\widehat{u}(x,t-h)}{h}.$$
Substituting the above test function $\varphi(x,t)$ into \eqref{limitu}, we obtain
\begin{align}\label{testlimit1}
	\int_{0}^{\infty}& \int_{\Omega}\xi(t) \gamma(\Theta)\nabla v \cdot(S_h\nabla\widehat{v})(x,t)\rd x\rd t\notag\\
	&=\int_{0}^{\infty} \int_{\Omega} v \varphi_t\rd x\rd t
	+ \int_{\Omega} u_{0t} \varphi(\cdot, 0)\rd x
	-a \int_{0}^{\infty} \int_{\Omega} \nabla u \nabla\varphi \rd x\rd t
	+ \int_{0}^{\infty} \int_{\Omega} f(\Theta) \nabla\varphi\rd x\rd t\notag\\
	&=\int_{0}^{\infty} \int_{\Omega} v\cdot\xi(t)\frac{\widehat v(x,t)-\widehat v(x,t-h)}{h}\rd x\rd t+ \int_{0}^{\infty} \int_{\Omega} v\cdot\xi'(t)(S_h\widehat{v})(x,t)\rd x\rd t+\int_{\Omega} u_{0t}^2\rd x\notag\\
	&~~~  -a \int_{0}^{\infty} \int_{\Omega} \nabla u \cdot\xi(t)(S_h\nabla \widehat{v})(x,t) \rd x\rd t
	+ \int_{0}^{\infty} \int_{\Omega} f(\Theta) \cdot\xi(t)(S_h\nabla\widehat{v})(x,t)\rd x\rd t.
\end{align}
Here we have used $\varphi(x,0)=u_{0t}$ since the fact that $\xi(0)=1$ and $(S_h\widehat{v})(x,0)=\frac1h\int_{-h}^0u_{0t}(x)\rd s=u_{0t}$ with \eqref{hatv} for a.e. $x\in\Omega$. 
Combining \eqref{hatv}, Cauchy-Schwarz inequality and the fact that
$$\frac{\xi(t+h)-\xi(t)}{h}\to \xi('t)\quad\text{and}\quad\frac1h\int_{0}^h\xi(t)\rd t\to\xi(0)=1,$$
as $h\searrow0$, then we get
\begin{align*}
	\int_{0}^{\infty} \int_{\Omega}\xi(t)v(x,t)&\frac{\widehat v(x,t)-\widehat v(x,t-h)}{h}\rd x\rd t\notag\\
	&=\frac1h\int_{0}^{\infty}\int_{\Omega}\xi(t) v^2(x,t)\rd x\rd t
	-\frac1h\int_{0}^{\infty}\int_{\Omega}\xi(t)v(x,t)\widehat{v}(x,t-h)\rd x\rd t\notag\\
	&\geqslant\frac{1}{2h}\int_{0}^{\infty}\int_{\Omega}\xi(t)v^2(x,t)\rd x\rd t
	-\frac{1}{2h}\int_{0}^{\infty}\int_{\Omega}\xi(t)\widehat{v}^2(x,t-h)\rd x\rd t\notag\\
	&=\frac{1}{2h}\int_{0}^{\infty}\int_{\Omega}\xi(t)v^2(x,t)\rd x\rd t-\frac{1}{2h}\int_{-h}^{\infty}\int_{\Omega}\xi(t+h)v^2(x,t)\rd x\rd t\notag\\
	&=-\frac12\int_0^\infty\int_\Omega\frac{\xi(t+h)-\xi(t)}{h}v^2(x,t)\rd x\rd t
	-\frac{1}{2h}\int_{0}^{h}\xi(t)\rd t\int_{\Omega}u_{0t}^2\rd x.
\end{align*}
Therefore we have 
$$\limsup_{h\searrow0} \int_{0}^{\infty} \int_{\Omega} v\cdot\xi(t)\frac{\widehat v(x,t)-\widehat v(x,t-h)}{h}\rd x\rd t
\geqslant-\frac12\int_0^\infty\int_\Omega\xi'(t)v^2\rd x \rd t-\frac12\int_{\Omega}u_{0t}^2\rd x.$$
 From \eqref{weak*vuTheta}, Lemma \ref{LemSteklov}, $\xi(t)\in C_0^\infty([0,\infty))$ and $\xi(0)=1$, together with the above inequality, returning to \eqref{testlimit1}, we arrive at
\begin{align}\label{nablauv}
	\int_{0}^{\infty}\int_{\Omega}\xi(t) \gamma(\Theta)|\nabla v|^2\rd x\rd t
	&\geqslant\frac12\int_0^\infty\int_\Omega\xi'(t)v^2\rd x \rd t+\frac12\int_{\Omega}u_{0t}^2\rd x+ \int_{0}^{\infty} \int_{\Omega}\xi(t)f(\Theta)\nabla v\rd x\rd t\notag\\
	&~~~  -a\limsup_{h\searrow0} \int_{0}^{\infty} \int_{\Omega}\xi(t)\nabla u(S_h\nabla\widehat{v})(x,t)\rd x\rd t.
\end{align}
Combining Lemma \ref{LemSteklov} with the Cauchy-Schwarz inequality again, it yields that
\begin{align*}
	\int_{0}^{\infty} \int_{\Omega}\xi(t)&\nabla u(S_h\nabla\widehat{v})(x,t)\rd x\rd t\\
	&=\int_{0}^{\infty} \int_{\Omega}\xi(t)\nabla  u\frac{\nabla\widehat{u}(x,t) - \nabla\widehat{u}(x,t-h)}{h}\rd x\rd t\\
	&=\frac1h\int_{0}^{\infty} \int_{\Omega} \xi(t)\left(\nabla\widehat{u}(x,t) - \nabla\widehat{u}(x,t-h)\right)^2\rd x\rd t\\
	&~~~+\frac1h\int_{0}^{\infty} \int_{\Omega} \xi(t)\nabla\widehat{u}(x,t-h)\left(\nabla\widehat{u}(x,t) - \nabla\widehat{u}(x,t-h)\right)\rd x\rd t\\
	&=h\int_{0}^{\infty} \int_{\Omega} \xi(t)\big(S_h\nabla
	\widehat{v}\big)^2\rd x\rd t-\frac1h\int_{0}^{\infty} \int_{\Omega}\xi(t)|\nabla\widehat{u}(x,t-h)|^2\\
	&~~~+\frac1h\int_{0}^{\infty} \int_{\Omega} \xi(t)\nabla\widehat{u}(x,t-h)\nabla\widehat{u}(x,t)\rd x\rd t\\
	&\leqslant h\int_{0}^{\infty} \int_{\Omega} \xi(t)\big(S_h\nabla\widehat{v}\big)^2\rd x\rd t-\frac{1}{2h}\int_{-h}^{\infty} \int_{\Omega}\xi(t+h)|\nabla\widehat{u}(x,t)|^2\rd x\rd t\\
	&~~~ +\frac{1}{2h}\int_{0}^{\infty} \int_{\Omega} \xi(t)|\nabla\widehat{u}(x,t)|^2\rd x\rd t\\
	&=h\int_{0}^{\infty} \int_{\Omega} \xi(t)\big(S_h\nabla\widehat{v}\big)^2\rd x\rd t
	-\frac{1}{2h}\int_{-h}^{0} \int_{\Omega}\xi(t+h)|\nabla \widehat u(x,t)|^2\rd x\rd  t\\
	&~~~-\frac{1}{2}\int_{0}^{\infty} \int_{\Omega} \frac{\xi(t+h)-\xi(t)}{h}|\nabla u(x,t)|^2\rd x\rd t.
\end{align*}
Then recalling Lemma \ref{LemSteklov} again, we obtain
\begin{align}\label{supuv}
	\limsup_{h\searrow0}\int_{0}^{\infty} \int_{\Omega}\xi(t)\nabla u(S_h\nabla\widehat{v})(x,t)\rd x\rd t
	\leqslant-\frac12\int_\Omega|\nabla u_0|^2\rd x-\frac12\int_0^\infty\int_\Omega\xi'(t)|\nabla u|^2\rd x\rd t.
\end{align}
Here we have used the fact that 
$$\int_{0}^{\infty} \int_{\Omega} \xi(t)\big(S_h\nabla\widehat{v}\big)^2\rd x\rd t<\infty,\quad\frac{\xi(t+h)-\xi(t)}{h}\to \xi('t),\quad \text{as}~h\searrow0,$$
and $$\frac1h\int_{-h}^0\xi(t+h)\rd t\to\xi(0)=1,\quad \text{as}~h\searrow0.$$
Substituting \eqref{supuv} into \eqref{nablauv}, we have
\begin{align}\label{I12*}
	\int_{0}^{\infty}\int_{\Omega}\xi(t) \gamma(\Theta)|\nabla v|^2\rd x\rd t&-\frac a2\int_{0}^{\infty} \int_{\Omega}\xi'(t)|\nabla u|^2\rd x\rd t\notag\\
	&\geqslant	\frac12\int_{\Omega}u_{0t}^2\rd x+\frac a2\int_\Omega|\nabla u_0|^2\rd x+\frac12\int_0^\infty\int_\Omega\xi'(t)v^2\rd x \rd t\notag\\
	&~~~ + \int_{0}^{\infty} \int_{\Omega}\xi(t)f(\Theta)\nabla v\rd x\rd t.
\end{align}

\textbf{Step 2.} 
Multiplying the first equation in \eqref{equ1*} by $\xi(t)v_\varepsilon$ and integrating over $\Omega\times(0,\infty)$ by parts, we get
\begin{align*}
	\frac{d}{dt} \Bigg\{ \frac{1}{2} \int_0^\infty\int_\Omega \xi(t)v_{\varepsilon}^2\rd x\rd t& + \frac{a}{2} \int_0^\infty\int_\Omega \xi(t) |\nabla u_{\varepsilon}|^2\rd x\rd t \Bigg\}-\frac{1}{2} \int_0^\infty\int_\Omega \xi'(t)v_{\varepsilon}^2\rd x\rd t\\
	&- \frac{a}{2} \int_0^\infty\int_\Omega \xi'(t)|\nabla u_{\varepsilon}|^2\rd x\rd t+\int_0^\infty\int_{\Omega}\xi(t) \gamma_\varepsilon(\Theta_\varepsilon)|\nabla v_{\varepsilon}|^2\rd x\rd t\\
	&+\varepsilon\int_0^\infty\int_\Omega \xi(t)|\Delta v_\varepsilon|^2\rd x\rd t+\varepsilon a\int_0^\infty\int_\Omega\xi(t)|\Delta u_{\varepsilon}|^2\rd x\rd t\\
	&=\int_0^\infty\int_{\Omega}\xi(t) f_\varepsilon(\Theta_\varepsilon)\nabla v_{\varepsilon}\rd x\rd t.
\end{align*}
Dropping some negative terms on the left side of the above inequality, we have
\begin{align*}
	\int_0^\infty\int_{\Omega}\xi(t) \gamma_\varepsilon(\Theta_\varepsilon)|\nabla v_{\varepsilon}|^2\rd x\rd t&-\frac{a}{2} \int_0^\infty\int_\Omega \xi'(t)|\nabla u_{\varepsilon}|^2\rd x\rd t\notag\\
	&\leqslant\int_0^\infty\int_{\Omega}\xi(t) f_\varepsilon(\Theta_\varepsilon)\nabla v_{\varepsilon}\rd x\rd t+\frac{1}{2} \int_0^\infty\int_\Omega \xi'(t)v_{\varepsilon}^2\rd x\rd t\notag\\
	&~~~ +\frac12\int_\Omega |v_{0\varepsilon}|^2\rd x+\frac a2\int_\Omega|\nabla u_{0\varepsilon}|^2\rd x.
\end{align*}
From Lemma \ref{Lemwekav}, \eqref{smoothdata} and $\xi(t),~\xi'(t)\in C_0^\infty([0,\infty))$, recalling the above inequality, we obtain
\begin{align*}
	\limsup_{\varepsilon=\varepsilon_k\searrow0}\int_0^\infty\int_{\Omega}\xi(t) \gamma_\varepsilon(\Theta_\varepsilon)|\nabla v_{\varepsilon}|^2\rd x\rd t
	&-\frac{a}{2} \int_0^\infty\int_\Omega \xi'(t)|\nabla u_{\varepsilon}|^2\rd x\rd t\\
	&\leqslant\int_0^\infty\int_{\Omega}\xi(t) f(\Theta)\nabla v\rd x\rd t+\frac{1}{2} \int_0^\infty\int_\Omega \xi'(t)v^2\rd x\rd t\\
	&~~~ +\frac12\int_\Omega u_{0t}^2\rd x+\frac a2\int_\Omega|\nabla u_0|^2\rd x.
\end{align*}
Together with \eqref{I12*}, we get
\begin{align}\label{I12geq}
	\limsup_{\varepsilon=\varepsilon_k\searrow0}\bigg(\int_0^\infty\int_{\Omega}\xi(t) &\gamma_\varepsilon(\Theta_\varepsilon)|\nabla v_{\varepsilon}|^2\rd x\rd t-\frac{a}{2} \int_0^\infty\int_\Omega \xi'(t)|\nabla u_{\varepsilon}|^2\rd x\rd t\bigg)\notag\\
	&\leqslant \int_{0}^{\infty}\int_{\Omega}\xi(t) \gamma(\Theta)|\nabla v|^2\rd x\rd t-\frac12\int_0^\infty\int_\Omega\xi'(t)v^2\rd x \rd t.	
\end{align}

\textbf{Step 3.} 
Let us introduce the following notations
$$I_{1,\varepsilon}:=\int_{0}^{\infty}\int_{\Omega}\xi(t) \gamma_\varepsilon(\Theta_\varepsilon)|\nabla v_\varepsilon|^2\rd x\rd t,$$
and
$$I_{2,\varepsilon}:=-\frac{a}{2} \int_0^\infty\int_\Omega \xi'(t)|\nabla u_{\varepsilon}|^2\rd x\rd t.$$
Recalling weakly lower semicontinuity of $I_{1,\varepsilon}$ and $I_{2,\varepsilon}$ in $L^2\big((0,\infty)\times\Omega\big)$ and $\xi'(t)\leqslant0$ for all $t\geqslant0$, we get
\begin{equation}\label{I1lower}
	I_1:=\int_{0}^{\infty}\int_{\Omega}\xi(t) \gamma(\Theta)|\nabla v|^2\rd x\rd t\leqslant \liminf_{\varepsilon=\varepsilon_k\searrow0}I_{1,\varepsilon},
\end{equation}
and 
\begin{equation}\label{I2lower}
  I_2:=-\frac12\int_0^\infty\int_\Omega\xi'(t)|\nabla u|^2\rd x \rd t\leqslant \liminf_{\varepsilon=\varepsilon_k\searrow0}I_{2,\varepsilon}.
\end{equation}
Combining \eqref{I12geq}, \eqref{I1lower} and \eqref{I2lower}, it follows that
\begin{equation}\label{I12equ}
	\limsup_{\varepsilon=\varepsilon_k\searrow0}\left(I_{1,\varepsilon}+I_{2,\varepsilon}\right)\leqslant I_1+I_2\leqslant\liminf_{\varepsilon=\varepsilon_k\searrow0}\left(I_{1,\varepsilon}+I_{2,\varepsilon}\right).
\end{equation}
We claim $I_{1,\varepsilon}\rightarrow I_1$. Otherwise one can find a subsequence $\{\varepsilon_{k_j}\}_{j\in\mathbb{N}}$ of $\{\varepsilon_{k}\}_{k\in\mathbb{N}}\subset(0,1)$ and $c_1>0$ such that
$$I_{1,\varepsilon_{k_j}}\geqslant I_1+c_1.$$
Then, combining \eqref{I12geq} and \eqref{I12equ}, it yields
$$I_{2,\varepsilon_{k_j}}\leqslant I_2-c_1+o(1).$$
That is 
$$I_2\geqslant \liminf_{\varepsilon_{k_j}\searrow0} I_{2,\varepsilon_{k_j}}+c_1,$$
thereby contradicting \eqref{I2lower}.
Then, for any $T>0$, we get
$$\sqrt{\gamma_\varepsilon(\Theta_\varepsilon) }\nabla v_{\varepsilon}
	\to \sqrt{\gamma(\Theta) }\nabla v\quad\text{in}~L^2\left(\Omega\times(0,T)\right),$$
as $\varepsilon=\varepsilon_{k}\searrow0,~k\rightarrow\infty$, due to $\xi\equiv1$ in $(0,T)$. Thus the proof of Lemma \ref{Lemgammato} is finished now.
\end{proof}

We are now in a position to complete the proof of Theorem \ref{Th1}. Combining the convergence properties established in Lemmas \ref{Lemwekav}, \ref{LemSteklov} and \ref{Lemgammato}, we shall pass to the limit in the regularized system \eqref{equ1*} and verify that the limit functions $(v,u,\Theta)$ satisfy the weak formulation of the original problem \eqref{equ1} in the sense of Definition \ref{solution}.

\noindent\textbf{Proof of Theorem \ref{Th1}.} 
We take the functions $(v, u, \Theta)$ as constructed in Lemma \ref{Lemwekav}. The regularity property \eqref{regular} follows immediately from the statement \eqref{weak*vuTheta} in Lemma \ref{Lemwekav}.
The weak identity \eqref{equmoment} is equivalent to \eqref{limitu} by virtue of Lemma \ref{Lemwekav} and the relation $u_t = v$ from \eqref{vuta.e.}.
	
To complete the proof, it remains to verify that the weak formulation of the temperature equation \eqref{equentropy} is satisfied. To this end, we fix an arbitrary test function $\varphi \in C^{\infty}_{0}(\Omega\times [0,\infty))$. Multiplying the third equation of \eqref{equ1*} by the above test function $\varphi$ and integrating over $\Omega\times(0,T) (T>0)$ by parts, for all $\varepsilon\in(0,1)$, we have
\begin{align*}
	- \int_{0}^{\infty} \int_{\Omega} &\Theta_{\varepsilon} \varphi_{t}\rd x\rd t - \int_{\Omega} \Theta_{0\varepsilon} \varphi(\cdot, 0)\rd x\rd t\\
	& = - \int_{0}^{\infty} \int_{\Omega} \nabla\Theta_{\varepsilon}\nabla \varphi\rd x\rd t 
	+ \int_{0}^{\infty} \int_{\Omega} \gamma_{\varepsilon}(\Theta_{\varepsilon}) |\nabla v_{\varepsilon}|^{2} \varphi \rd x\rd t- \int_{0}^{\infty} \int_{\Omega} f_{\varepsilon}(\Theta_{\varepsilon}) \nabla v_{\varepsilon} \varphi\rd x\rd t.
\end{align*}
We now pass to the limit along the subsequence $\{\varepsilon_{k}\}_{k\in \mathbb{N}}$ provided by Lemma \ref{Lemwekav} and Lemma \ref{Lemgammato}. 
By the convergence $\Theta_{\varepsilon} \to \Theta$ in $L^{q}_{loc}(\Omega\times [0,\infty))$ for all $q \in \big[1,\frac{N+2}{N}\big)$ from Lemma \ref{Lemwekav} and the initial data convergence in \eqref{smoothdata}, we obtain
$$\int_{0}^{\infty} \int_{\Omega} \Theta_{\varepsilon} \varphi_{t}\rd x\rd t \to \int_{0}^{\infty} \int_{\Omega} \Theta \varphi_{t}\rd x\rd t, \quad \text{and} \quad \int_{\Omega} \Theta_{0\varepsilon} \varphi(\cdot, 0)\rd x\rd t \to\int_{\Omega} \Theta_{0} \varphi(\cdot, 0)\rd x\rd t.$$
The convergence $\nabla\Theta_{\varepsilon} \rightharpoonup \nabla\Theta$ in $L^{r}_{loc}(\Omega \times [0,\infty))$ for all $r \in \left(1, \frac{N+2}{N+1}\right)$ from \eqref{weak*vuTheta} yields
$$\int_{0}^{\infty} \int_{\Omega} \nabla\Theta_{\varepsilon} \nabla\varphi\rd x\rd t \to  \int_{0}^{\infty} \int_{\Omega} \nabla\Theta\nabla \varphi\rd x\rd t.$$
By recalling \eqref{vuta.e.} and the strong convergence result of Lemma \ref{Lemgammato}, we have directly (by taking the square and using the uniform boundedness of $\gamma$)
$$\int_{0}^{\infty} \int_{\Omega} \gamma_{\varepsilon}(\Theta_{\varepsilon}) |\nabla v_{\varepsilon}|^{2} \varphi\rd x\rd t \to \int_{0}^{\infty} \int_{\Omega} \gamma(\Theta) |\nabla u_t|^{2} \varphi\rd x\rd t.$$
According to Lemma \ref{Lemwekav} and \eqref{vuta.e.} again, we conclude
$$\int_{0}^{\infty} \int_{\Omega} f_{\varepsilon}(\Theta_{\varepsilon}) \nabla v_{\varepsilon} \varphi \rd x\rd t\to   \int_{0}^{\infty} \int_{\Omega} f(\Theta) \nabla u_t \varphi\rd x\rd t.$$
Collecting all the above limit relations, we deduce that the limiting function $(u, \Theta)$ satisfies the weak identity \eqref{equentropy} for every $\varphi \in C^{\infty}_{0}(\Omega \times [0,\infty))$. This completes the  proof that $(u, \Theta)$ is a global weak solution in the sense of Definition \ref{solution}.
$\hfill\square$


\begin{thebibliography}{99}
\bibitem{Amann1993}\label{Amann1993}
H. Amann, Nonhomogeneous linear and quasilinear elliptic and parabolic boundary value problems. In: Function spaces, differential operators and nonlinear analysis (Friedrichroda, 1992), Volume 133 of Teubner-Texte Math., Teubner, Stuttgart, 9--126 (1993).

\bibitem{BDGO}\label{BDGO}
D. Blanchard, O. Guib\'e, Existence of a solution for a nonlinear system in thermoviscoelasticity, Adv. Differ. Equ. {\bf 5}, 1221--1252 (2000). 

\bibitem{CH1994}\label{CH1994}
Z.~M. Chen, K.-H. Hoffmann, On a one-dimensional nonlinear thermoviscoelastic model for structural phase transitions in shape memory alloys, J. Differential Equations {\bf 112}(2), 325--350 (1994). 

\bibitem{WinklerMMM}\label{WinklerMMM}
L. Claes, J. Lankeit, M. Winkler, A model for heat generation by acoustic waves in piezoelectric materials: global large-data solutions, Math. Models Methods Appl. Sci. {\bf 35}(11), 2465--2512 (2025).

\bibitem{CWPrepint}\label{CWPrepint}
L. Claes, M. Winkler, Describing smooth small-data solutions to a quasilinear hyperbolic-parabolic system by $W^{1,p}$energy analysis, Nonlinear Anal. Real World Appl. {\bf 91}, 104580 (2026). 

\bibitem{Dafermos1982}\label{Dafermos1982}
C.~M. Dafermos, Global smooth solutions to the initial-boundary value problem for the equations of one-dimensional nonlinear thermoviscoelasticity, SIAM J. Math. Anal. {\bf 13}(3), 397--408 (1982).

\bibitem{DH1982}\label{DH1982}
C.~M. Dafermos, L. Hsiao, Global smooth thermomechanical processes in one-dimensional nonlinear thermoviscoelasticity, Nonlinear Anal. {\bf 6}(5), 435--454 (1982).

\bibitem{Steklov1993}\label{Steklov1993}
E. DiBenedetto, {\it Degenerate parabolic equations}, Universitext, Springer, New York, 1993.

\bibitem{Fricke2025}\label{Fricke2025}
T.~J. Fricke, Local and global solvability in a viscous wave equation involving general temperature-dependence, Acta Appl. Math. {\bf 200}(1),  (2025). 

\bibitem{Friesen2024}\label{Friesen2024}
O. Friesen, L. Claes, C. Scheidemann, N. Feldmann, T. Hemsel, B. Henning, Estimation of temperature-dependent piezoelectric material parameters using ring-shaped specimens, In: 2023 International Congress on Ultrasonics, Beijing, China, Vol. 2022, 012125. IOP Publishing, 2024.


\bibitem{GJZW}\label{GJZW}
J.~A.~P. Gawinecki, W.~M. Zajaczkowski, Global regular solutions to two-dimensional thermoviscoelasticity, Commun. Pure Appl. Anal. {\bf 15}(3), 1009--1028 (2016). 

\bibitem{GBSD2007}\label{GBSD2007}
Z. Gubinyi, C. Batur, A. Sayir, F. Dynys,  Electrical properties of PZT piezoelectric ceramic at high temperatures, J. Electroceram. {\bf 20}(2), 95--105 (2008).

\bibitem{GZ1999}\label{GZ1999}
B. Guo, P. Zhu, Global existence of smooth solution to nonlinear thermoviscoelastic system with clamped boundary conditions in solid-like materials, Commun. Math. Phys. {\bf 203}(2), 365--383 (1999).

\bibitem{Winkler2005JDE}\label{Winkler2005JDE}
D. Horstmann, M. Winkler, Boundedness vs. blow-up in a chemotaxis system, J. Differential Equations {\bf 215}(1), 52--107 (2005). 

\bibitem{HL1998}\label{HL1998}
L. Hsiao, T. Luo, Large-time behavior of solutions to the equations of one-dimensional nonlinear thermoviscoelasticity, Quart. Appl. Math. {\bf 56}(2), 201--219 (1998).  

\bibitem{JS1993}\label{JS1993}
S. Jiang, Global large solutions to initial-boundary value problems in one-dimensional nonlinear thermoviscoelasticity, Quart. Appl. Math. {\bf 51}(4), 731--744 (1993).

\bibitem{Kim1983}\label{Kim1983}
J.~U. Kim, Global existence of solutions of the equations of one-dimensional thermoviscoelasticity with initial data in ${\rm BV}$\ and $L\sp{1}$, Ann. Scuola Norm. Sup. Pisa Cl. Sci. (4) {\bf 10}(3), 357--427 (1983).

\bibitem{MART}\label{MART}
A. Mielke, T. Roub\'i\v cek, Thermoviscoelasticity in Kelvin-Voigt rheology at large strains, Arch. Ration. Mech. Anal. {\bf 238}(1), 1--45 (2020).

\bibitem{PIZW}\label{PIZW}
I. Paw\l ow, W.~M. Zajaczkowski, Global regular solutions to three-dimensional thermo-visco-elasticity with nonlinear temperature-dependent specific heat, Commun. Pure Appl. Anal. {\bf 16}(4), 1331--1371 (2017).

\bibitem{RZ1997}\label{RZ1997}
R. Racke, S. Zheng, Global existence and asymptotic behavior in nonlinear thermoviscoelasticity, J. Differential Equations {\bf 134}(1), 46--67 (1997).

\bibitem{RT2009}\label{RT2009}
T. Roub\'i\v cek, Thermo-visco-elasticity at small strains with $L^1$-data, Quart. Appl. Math. {\bf 67}(1),  47--71 (2009).

\bibitem{RT2013}\label{RT2013}
T. Roub\'i\v cek, Nonlinearly coupled thermo-visco-elasticity, NoDEA Nonlinear Differ. Equ. Appl. {\bf 20}(3), 1243--1275 (2013).

\bibitem{SY1995}\label{SY1995}
Y. Shibata, Global in time existence of small solutions of nonlinear thermoviscoelastic equations, Math. Methods Appl. Sci. {\bf 18}(11), 871--895 (1995).

\bibitem{SZZ1999}\label{SZZ1999}
W. Shen, S. Zheng, P. Zhu, Global existence and asymptotic behavior of weak solutions to nonlinear thermoviscoelastic systems with clamped boundary conditions, Quart. Appl. Math. {\bf 57}(1), 93--116 (1999). 

\bibitem{Watson2000}\label{Watson2000}
S.~J. Watson, Unique global solvability for initial-boundary value problems in one-dimensional nonlinear thermoviscoelasticity, Arch. Ration. Mech. Anal. {\bf 153}(1), 1--37 (2000).

\bibitem{WinklerAMO}\label{WinklerAMO}
M. Winkler, Rough data in an evolution system generalizing 1D thermoviscoelasticity with temperature-dependent parameters, Appl. Math. Optim. {\bf 91}(2), (2025). 

\bibitem{WinklerZAMP}\label{WinklerZAMP}
M. Winkler, Large-data solutions in one-dimensional thermoviscoelasticity involving temperature-dependent viscosities, Z. Angew. Math. Phys. {\bf 76}(5), (2025). 

\bibitem{ZS1989}\label{ZS1989}
S. Zheng, W. Shen, Global smooth solutions to Cauchy problem of equations of one-dimensional thermoviscoelasticity, J. Partial Differential Equations {\bf 2}(2), 26--38 (1989). 

\end{thebibliography}
\end{document}